\documentclass[12pt,oneside,reqno]{amsart}

\usepackage{amsfonts}
\usepackage{amsmath}
\usepackage{amssymb}
\usepackage{amsthm}
\usepackage{enumerate}
\usepackage{graphicx}
\usepackage{hyperref}
\usepackage{mathrsfs}
\usepackage{stmaryrd}
\usepackage{color}

\newcommand{\be}{\begin{eqnarray}}
\newcommand{\ee}{\end{eqnarray}}
\newcommand{\ce}{\begin{eqnarray*}}
\newcommand{\de}{\end{eqnarray*}}

\newtheorem{theorem}{Theorem}[section]
\newtheorem{proposition}[theorem]{Proposition}
\newtheorem{lemma}[theorem]{Lemma}
\newtheorem{corollary}[theorem]{Corollary}
\newtheorem{remark}[theorem]{Remark}
\newtheorem{definition}[theorem]{Definition}
\newtheorem{examples}[theorem]{Examples}
\newtheorem{assumption}[theorem]{Assumption}

\def\bt{\begin{theorem}}
\def\et{\end{theorem}}
\def\bp{\begin{proposition}}
\def\ep{\end{proposition}}
\def\bl{\begin{lemma}}
\def\el{\end{lemma}}
\def\bc{\begin{corollary}}
\def\ec{\end{corollary}}
\def\bd{\begin{definition}}
\def\ed{\end{definition}}
\def\br{\begin{remark}}
\def\er{\end{remark}}
\def\bx{\begin{examples}}
\def\ex{\end{examples}}
\def\ba{\begin{assumption}}
\def\ea{\end{assumption}}

\def\e{\varepsilon}

\def\[{{\Big[}}
\def\]{{\Big]}}
\def\<{{\langle}}
\def\>{{\rangle}}
\def\({{\Big(}}
\def\){{\Big)}}

\def\geq{\geqslant}
\def\leq{\leqslant}

\def\no{\nonumber}

\def\cB{{\mathcal B}}
\def\cC{{\mathcal C}}
\def\cD{{\mathcal D}}
\def\cE{{\mathcal E}}

\def\cG{{\mathcal G}}
\def\cH{{\mathcal H}}

\def\cL{{\mathcal L}}
\def\cM{{\mathcal M}}

\def\cO{{\mathcal O}}
\def\cP{{\mathcal P}}

\def\cR{{\mathcal R}}
\def\cS{{\mathcal S}}

\def\cX{{\mathcal X}}

\def\mD{{\mathbb D}}
\def\mE{{\mathbb E}}

\def\mN{{\mathbb N}}

\def\mP{{\mathbb P}}

\def\mR{{\mathbb R}}

\allowdisplaybreaks

\begin{document}

\title{Bismut-Elworthy-Li Formulae for Forward-Backward SDEs with Jumps and Applications}

\author{Jiagang Ren$^{1}$, Hua Zhang$^{2}$}


\date{}

\dedicatory{$^{1}$School of Mathematics, Sun Yat-Sen University,\\
Guangzhou, Guangdong 510275, P.R.China\\
$^{2}$School of Statistics and Data Science, Jiangxi University of Finance and Economics,\\
Nanchang, Jiangxi 330013, P.R.China\\
Emails: J. Ren: renjg@mail.sysu.edu.cn\\
H. Zhang: zh860801@163.com}


\thanks{This work is supported by National Natural Science Foundation of China (Grant Nos. 12261038 and 12371152), and Natural Science Foundation of Jiangxi Province (Grant Nos. 20232BAB201004 and 20242BAB23003).}

\begin{abstract}
Under nondegeneracy assumptions on the diffusion coefficients, we establish the derivative formulae of Bismut-Elworthy-Li's type for forward-backward stochastic differential equations with respect to Poisson random measure using the lent particle method created by Bouleau and Denis, which is not given before. Applying this formula, the existence and uniqueness of a solution of nonlocal quasi-linear integral partial differential equations, which are differentiable with respect to the space variable, are obtained, even if the initial datum and coefficients of this equation are not.
\end{abstract}

\maketitle

\section{Introduction}

In the present paper, we are concerned with the problem of derivative formulae of Bismut-Elworthy-Li's type for the following forward-backward stochastic differential equations (SDEs in short) with respect to Poisson random measure:
\be\label{SDE}
\begin{cases}
X(\tau,t,x)=x+\int_t^{\tau}b(s,X(s,t,x))ds\\
\quad\quad\quad\quad\quad+\int_t^{\tau}\int_{\cO}\sigma(s,X(s-,t,x))u\widetilde{N}(ds,du),\tau\in[t,T],\\
X(\tau,t,x)=x,\tau\in[0,t],
\end{cases}
\ee
where $b:[0,T]\times\mR^d\rightarrow\mR^d$, $\sigma:[0,T]\times\mR^d\rightarrow\mR^d\times\mR^l$,
$\cO=\{u\in\mR^l;|u|\leq1\}\setminus\{0\}$, $\widetilde{N}$ is the compensated process of some Poisson point process $N$ on $(\cO,\cB(\cO))$ with intensity $\nu$, and
\be\label{BSDE}
\begin{split}
&Y(\tau,t,x)+\int_{\tau}^T\int_{\cO}Z(s-,t,x,u)\widetilde{N}(ds,du)\\
=&-\int_{\tau}^T\psi(s,X(s,t,x),Y(s,t,x),\int_{\cO}Z(s,t,x,u)l(u)\nu(du))ds+\phi(X(T,t,x)),
\end{split}
\ee
where $\tau\in[0,T]$, $\psi:[0,T]\times\mR^d\times\mR\times\mR\rightarrow\mR$ and $\phi:\mR^d\rightarrow\mR$ are given functions, and $l:\cO\rightarrow\mR$ is a Borel function with growth condition $|l(u)|\leq C(1\wedge|u|)$, $u\in\cO$.

A derivative formula for Wiener functionals was first established by Bismut in \cite{B} using Malliavin calculus on Wiener space. This result was subsequently extended by Elworthy and Li in their foundational work \cite{EL}, where they systematically derived formulae for the derivatives of heat semigroups—this contribution, together with Bismut’s earlier result, laid the core theoretical basis for what is now known as the Bismut-Elworthy-Li formulae. Unlike the Wiener space, there are several versions of Malliavin calculus on Poisson space which are not equivalent. In the recent years, quite a few papers appeared concerning the Bismut-Elworthy-Li formulae for SDEs with jumps, using various approaches, e.g., SDEs driven by general Poisson jumps processes by using stochastic diffeomorphism flows and Girsanov's transformation (\cite{T}), SDEs driven by subordinated Brownian motion by using conditional Malliavin calculus (\cite{Z,WXZ2}), SDEs driven by $\alpha$-stable-like noise by using Bismut's approach(\cite{WXZ1}).

In a series of works, Bouleau and Denis (see \cite{BD3} and references therein) presented systematically a new method called the lent particle one which may be considered as a differential calculus on Poisson space completely comparable to the Malliavin differential calculus on the Wiener space. The main feature of this method is its simplicity and convenience when applied to Poisson functionals, especially to solutions of stochastic differential equations driven by L\'{e}vy processes. Using this method, the Bismut-Elworthy-Li formulae for SDEs with jumps is also established in our recent paper \cite{RZ} which generalizes and improves some previous works in this respect.

Although some authors have studied the Bismut-Elworthy-Li formulae for SDEs with jumps, no one has yet discussed the same problem for forward-backward SDEs with jumps. For the cases of forward-backward SDEs without jumps, under nondegeneracy assumptions on the diffusion coefficient (that may be nonconstant), Fuhrman and Tessitore \cite{FT} have proved an analogue of the Bismut-Elworthy-Li formulae. After that, Masiero \cite{M} has also obtained a Bismut-Elworthy-Li formula for forward-backward SDEs in a Markovian framework when the generator has quadratic growth with respect to $Z$. In the present paper, using the lent particle method, we establish a derivative formula of Bismut-Elworthy-Li's type for forward-backward SDEs with respect to Poisson random measure under the nondegeneracy assumption on the diffusion coefficients, which has not been given so far to our knowledge, and we apply this formula to the study of nonlocal quasi-linear integral partial differential equations (PDEs in short).

The organization of this paper is as follows. In Section 2, we recall the lent particle method which is the basic tool in our article. In Section 3, some notations and preliminary results on forward-backward SDEs are introduced. In Section 4, we obtain a derivative formula of Bismut-Elworthy-Li's type for forward-backward SDEs with jumps. Building upon this formula, in Section 5, we establish the existence and uniqueness of solutions to nonlocal quasi-linear integral-PDEs that are differentiable with respect to the space variable, regardless of the potential non-differentiability of the initial conditions or coefficients of the equation.

In the present paper, let $C_p$ be a positive constant only depending on some parameter $p$, whose value may change from line to line. When we don't want to emphasize this dependence we just use $C$ instead.

\section{Lent Particle Method}

\subsection{Set-up}

Let us first specify the general set-up in which we will work. We follow \cite{BD3} to which we refer for more details.

\subsubsection{Dirichlet Structure on the Bottom Space}

We start from a bottom space $(\Xi,\cG,\nu)$, where $\Xi$ is a separable Hausdorff space, $\cG$ its Borel $\sigma$-algebra and $\nu$ a $\sigma$-finite and diffuse measure on $(\Xi,\cG)$. Let $(\mathbf{d},e)$ be a local symmetric Dirichlet form on $L^2(\nu)$ which admits a carr\'{e} du champ operator $\gamma$. That is to say, $\gamma$ is the unique positive, symmetric and continuous bilinear form from $\mathbf{d}\times\mathbf{d}$ to $L^1(\nu)$ such that
\ce
e(f,g)=\frac12\int_{\Xi}\gamma[f,g]d\nu,\quad\forall f,g\in\mathbf{d}.
\de
Moreover, we suppose that there exist $\{k_n,n\in\mN\}\subset\mathbf{d}$ and $A_n\uparrow\cX$ such that $k_n1_{A_n}\uparrow 1$ and $\gamma[k_n]1_{A_n}=0$. The structure $(\Xi,\cG,\nu,\mathbf{d},\gamma)$ is called the bottom structure.

Since $\mathbf{d}$ is separable, the bottom Dirichlet structure admits a gradient operator, i.e., there exists a separable Hilbert space $H$ and a linear map $D$ from $\mathbf{d}$ into $L^2(\nu;H)$ such that
\ce
\gamma[u]=\|Du\|_H^2,\quad\forall u\in\mathbf{d}.
\de
Let $(R,\cR,\rho)$ be another probability space such that the vector space $L^2(R,\cR,\rho)$ is infinite dimensional. Take $H=L_0^2(R,\cR,\rho)=\{g\in L^2(R,\cR,\rho);\int_Rg(r)\rho(dr)=0\}$. The corresponding gradient will be denoted by $\flat$, and we assume without any loss of generality that constants belong to $\mathbf{d}_{loc}$ (see \cite[Chapter I, Definition 7.1.3]{BF}) and so that $1^{\flat}=0$.

\subsubsection{Space-time Setting and Dirichlet Structure on the Upper Space}

From now on we set $X=[0,T]\times \Xi$, $\cX=\cB([0,T])\times \cG$ and $\mu=dt\times\nu$. Define the Dirichlet structure on $(X,\cX,\mu)$ to be the product of the trivial one on $(L^2([0,T],dt),0)$ and $(\mathbf{d}, e)$ and we keep the same notations $\mathbf{d}$, $e$, $\flat$, $\gamma$ and $\mathbf{a}$, etc, for operators corresponding to this new Dirichlet form but note that they act only on the second variable.

It is known that $X$ is totally ordered (see \cite[Section 1, Theorem 11]{DM}) and we denote by $\prec$ such a total order relation. Set
\ce
\Omega:=\{\omega=\sum_{i=1}^\infty\e_{y_i};\quad y_i\in X,\quad\forall i,\quad\text{and}\quad y_1\prec y_2\prec\cdots\prec y_n\prec\cdots\}.
\de
Let $N$ be the Poisson random measure with intensity $\mu$ defined on $(\Omega,\mathcal{F},\mP)$ where $N(\omega)=\omega$, $\mathcal{F}$ is the $\sigma$-algebra generated by $N$ and $\mP$ the law of $N$. Define $\{\mathcal{F}_t,t\in[0,T]\}$ to be the $\mP$-complete right continuous filtration generated by $\{N([0,t]\times G), t\in[0,T], G\in\cG\}$.

We now introduce the creation and annihilation operator $\varepsilon^+$ and $\varepsilon^-$:
\ce
\forall(t,y)\in\mR^+\times\Xi,\quad\forall\omega\in\Omega,\\
\varepsilon_{(t,y)}^+(\omega)=\omega1_{\{(t,y)\in supp\omega\}}+(\omega+\varepsilon_{(t,y)})1_{\{(t,y)\notin supp\omega\}},\\
\forall(t,y)\in\mR^+\times\Xi,\quad\forall\omega\in\Omega,\\
\varepsilon_{(t,y)}^-(\omega)=\omega1_{\{(t,y)\notin supp\omega\}}+(\omega-\varepsilon_{(t,y)})1_{\{(t,y)\in supp\omega\}}.
\de
Denote $\mP_N:=\mP(d\omega)N_{\omega}(dt,dy)$.  Then it is well known (see \cite[Lemma 4.2]{BD3}) that the map $(\omega,(t,y))\mapsto (\varepsilon_{(t,y)}^+\omega,(t,y))$
sends $\mP_N$-negligible sets to $\mP\times\mu$-negligible ones, and the map $(\omega,(t,y))\mapsto (\varepsilon_{(t,y)}^-\omega,(t,y))$
sends $\mP\times\mu$-negligible sets to $\mP_N$-negligible ones.

If $N(\omega)=\sum_{i=1}^{\infty}\varepsilon_{y_i}$, then define
\ce
N\odot\rho(\omega,\hat{\omega}):=\sum_{i=1}^{\infty}\varepsilon_{(y_i,r_i(\hat{\omega}))},
\de
where $(r_i)$ is a sequence of i.i.d random variables independent of $N$ whose common law is $\rho$ and which are defined on some probability space $(\widehat{\Omega},\widehat{\mathcal{F}},\widehat{\mP})$. Hence $N\odot\rho$ is defined on the product probability pace $(\Omega,\mathcal{F},\mP)\times(\widehat{\Omega},\widehat{\mathcal{F}},\widehat{\mP})$. It is a Poisson random measure on $X\times R$ with compensator $\mu\times\rho$ which is called the marked Poisson random measure.

The marked Poisson random measure $N\odot\rho$ has the following useful features in connection with $N$ and the measure $\widehat{\mP}$ which is crucial in the subsequent sections. The proof is similar to that of \cite[Corollary 12]{BD1}, and hence we omit it.

\bp\label{marked Poisson measure}
Let $F$, $G:\Omega\times X\times R\rightarrow\mR$ be $\mathcal{F}\otimes\cX\otimes\cR$-measurable functions such that
\ce
&&\int_0^T\int_{\Xi}\int_R(|F(t,y,r)|+|F(t,y,r)|^2)\rho(dr)N(dt,dy)<\infty,\quad\mP-\text{a.e.},\\
&&\int_RF(t,y,r)\rho(dr)=0,\quad\mP_N-\text{a.e.}
\de
and
\ce
&&\int_0^T\int_{\Xi}\int_R(|G(t,y,r)|+|G(t,y,r)|^2)\rho(dr)N(dt,dy)<\infty,\quad\mP-\text{a.e.},\\
&&\int_RG(t,y,r)\rho(dr)=0,\quad\mP_N-\text{a.e.}
\de
Then the following relation holds $\mP_N$-a.e.
\ce
&&\widehat{\mE}[(\int_0^T\int_{\Xi}\int_RF(t,y,r)N\odot\rho(dt,dy,dr))(\int_0^T\int_{\Xi}\int_RG(t,y,r)N\odot\rho(dt,d\nu,dr))]\\
&=&\int_0^T\int_{\Xi}\int_RF(t,y,r)G(t,y,r)\rho(dr)N(dt,dy).
\de
Here and in the sequel, $\int_{\alpha}^{\beta}=\int_{(\alpha,\beta]}$ for $\alpha<\beta$, and $\widehat{\mE}$
denotes the expectation with respect to $\widehat{\mP}$.
\ep

On the product structure $(X_n,\cX_n,\mu_n,\mathbf{d}_n,\gamma_n)
=(X,\cX,\mu,\mathbf{d},\gamma)^n$ the Dirichlet form is defined by
\ce
e_n(f)=\frac12\int_{X_n}\gamma_n[f]d\mu_n.
\de
We refer to \cite[Chapter V]{BF} for a detailed account of the theory of this product Dirichlet structure. Let $\{p_t;t\in\mR^+\}$ be the semigroup associated to $e$ and $\{p^n_t;t\in\mR^+\}$ that associated to $e_n$. For $F\in L^2(\mP)$ with the chaos decomposition
\ce
F=\mE[F]+\sum_{n=1}^{\infty}I_n(f_n),
\de
define
\ce
P_tF=\mE[F]+\sum_{n=1}^\infty I_n(p^n_t(f_n)).
\de
Then $\{P_t;t\in\mR^+\}$ on $L^2(\mP)$  forms a symmetric strongly continuous
semigroup whose infinitesimal generator is given by
\ce
\cD(A)=\{F\in L^2(\mP);\lim_{t\downarrow 0}\frac{P_tF-F}{t}\quad\text{exists in}\quad L^2(\mP)\}
\de
and
\ce
A[F]=\lim_{t\downarrow 0}\frac{P_tF-F}{t},\quad\forall F\in\cD(A).
\de
Define
\ce
\mD=\{F\in L^2(\mP);\lim_{t\downarrow 0}\<\frac{F-P_tF}{t},F\>_{L^2(\mP)}<\infty\}
\de
and
\ce
\cE(F)=\lim_{t\downarrow0}\<\frac{F-P_tF}{t},F\>_{L^2(\mP)},\quad\forall F\in\mD.
\de
Then $(\mD,\cE)$ is a local symmetric Dirichlet form on $L^2(\mP)$ with a carr\'e du champ operator $\Gamma$ given by
\ce
\Gamma[F,G]=\widehat{\mE}[F^{\sharp}G^{\sharp}],
\de
where
\be\label{gradient}
F^\sharp=\int_0^T\int_{\Xi}\int_R\varepsilon^-(\varepsilon^+F)^{\flat}dN\odot\rho.
\ee
Moreover, we have
\ce
\Gamma[F]=\widehat{\mE}[(F^{\sharp})^2]=\int_0^T\int_{\Xi}\varepsilon^-(\gamma[\varepsilon^+F])dN,\quad\forall F\in\mD.
\de

We now recall the definition of the divergence operator $\delta_{\sharp}$.

\bd
The operator $\delta_{\sharp}:L^2(\mP\times\widehat{\mP})\rightarrow\mD$ is defined as the adjoint operator of the gradient $F\in\mD\mapsto F^{\sharp}\in L^2(\mP\times\widehat{\mP})$, and we denote by $dom\delta_{\sharp}$ the domain of $\delta_{\sharp}$.
\ed

We will need the following result which is due to \cite[Proposition 5.6]{BD3}.

\bp\label{divergence}
Let $F=F(\omega,\hat{\omega})$ be $\mathcal{F}^{\odot}$-measurable, where $\mathcal{F}^{\odot}$ is the $\sigma$-field on $\Omega\times\widehat{\Omega}$ generated by $N\odot\rho$, $F\in dom\delta_{\sharp}$, then
\ce
\delta_{\sharp}F=\int_0^T\int_{\Xi}\varepsilon_{(\alpha,y)}^-(\widehat{\mE}[\delta_{\flat}(\varepsilon_{(\alpha,y,r)}^+F)(\alpha,y)])N(d\alpha,dy),
\de
where $\varepsilon^+$ is relative to the Poisson random measure $N\odot\rho$ under $\mP\times \widehat{\mP}$ hence adds a point $(\alpha, y,r)$ while operator $\varepsilon$ is relative to $N$.
\ep

\section{Notations and Preliminary Results on the Forward-Backward Systems}

\subsection{Notations}

Throughout this paper, $\mR$ and $\mR^+$ are the set of real numbers and non-negative real numbers respectively, and $\mN$ and $\mN^*$ are the set of integers and positive integers respectively. We use the notation $|\cdot|$ to denote the usual Euclidean norm in $\mR^d$ and $\mR^d\otimes\mR^l$ simultaneously. From this section on, we denote $(\Xi,\cG,\nu)=(\cO,\cB(\cO),kdu)$.

For any real and separable Hilbert space $K$, the symbol $\cS^p(K)$, $1\leq p<\infty$, or $\cS^p$ where no confusion is possible, denotes the space of all predictable processes $\{Y(\tau);\tau\in[0,T]\}$, with values in $K$, such that the norm
\ce
\|Y\|_{\cS^p}=\mE[\sup_{\tau\in[0,T]}|Y({\tau})|^p]^{1/p}
\de
is finite. $\cS^{\infty}(K)$, or $\cS^{\infty}$ where no confusion is possible, denotes the space of all bounded predictable processes. 
We also denote by $\cM^p(K)$ ($1\leq p<\infty$), or $\cM^p$ where no confusion is possible, the space of all predictable processes $\{Z(\tau,u);\tau\in[0,T],u\in\cO\}$ with values in $K$, normed by
\ce
\|Z\|_{\cM^p}=\mE[(\int_0^T\int_{\cO}|Z(\tau,u)|^2\nu(du)d\tau)^{p/2}]^{1/p}.
\de

Given two finite-dimensional spaces $\mathbb{R}^m$ and $\mathbb{R}^n$, we say that a mapping $F: \mathbb{R}^m \rightarrow \mathbb{R}^n$ belongs to the class $\mathcal{C}^1(\mathbb{R}^m; \mathbb{R}^n)$ if it is continuous, continuously differentiable on $\mathbb{R}^m$, and its derivative $\nabla F: \mathbb{R}^m \rightarrow L(\mathbb{R}^m; \mathbb{R}^n)$ is continuous. For another finite-dimensional space $\mathbb{R}^p$, a mapping $F: \mathbb{R}^m \times \mathbb{R}^p \rightarrow \mathbb{R}^n$ belongs to the class $\mathcal{C}^{0,1}(\mathbb{R}^m \times \mathbb{R}^p; \mathbb{R}^n)$ if it is continuous, continuously differentiable with respect to $y \in \mathbb{R}^n$, and the partial derivative $\nabla_y F: \mathbb{R}^m \times \mathbb{R}^p \rightarrow L(\mathbb{R}^p; \mathbb{R}^n)$ is continuous. In this work, we focus on specific instances of these classes, such as $\mathcal{C}^{0,1}([0,T] \times \mathbb{R}^d; \mathbb{R})$ (functions continuous in time and Lipschitz in space), $\mathcal{C}^{1,1,1}(\mathbb{R}^d \times \mathbb{R} \times \mathbb{R}; \mathbb{R})$ (functions continuously differentiable in all three variables), and $\mathcal{C}^1(\mathbb{R}^d; \mathbb{R})$ (continuously differentiable functions).

We denote by $\cH_{\mD}$ the set of real-valued processes $\{Y(t);t\in[0,T]\}$ which satisfy
\ce
Y(t)\in\mD\quad\text{and}\quad\|Y\|_{\cH_{\mD}}:=\mE[\int_0^T|Y(t)|^2dt]+\mE\widehat{\mE}[\int_0^T|Y^{\sharp}(t)|^2dt]<\infty.
\de
We denote by $\cH_{\mD,\nu}$ the set of real-valued processes $\{Z(t,u);t\in[0,T],u\in\cO\}$ which satisfy $Z(t,u)\in\mD$ and
\ce
\|Z\|_{\cH_{\mD,\nu}}:=\mE[\int_0^T\int_{\cO}|Z(t,u)|^2\nu(du)dt]+\mE\widehat{\mE}[\int_0^T\int_{\cO}|Z^{\sharp}(t,u)|^2\nu(du)dt]<\infty.
\de
We also denote by $\cH_{\mD}^d$ and $\cH_{\mD,\nu}^d$ the space of $\mR^d$-valued processes such that each coordinate belongs to $\cH_{\mD}$ and $\cH_{\mD,\nu}$ respectively and we equip it with the standard norm of product topology.

\subsection{Forward SDEs}

The coefficients $b$ and $\sigma$ are assumed to satisfy the following conditions.

\begin{enumerate}
\item[\textbf{(FL)}]
There exists a constant $C_{FL}>0$ such that
\ce
|b(t,x)|\vee|\sigma(t,x)|\leq C_{FL}(1+|x|),\quad\forall t\in[0,T],~\forall x\in\mR^d
\de
and
\ce
|b(t,x)-b(t,y)|\vee|\sigma(t,x)-\sigma(t,y)|\leq C_{FL}|x-y|,\quad\forall t\in[0,T],~\forall x,y\in\mR^d.
\de
\item[\textbf{(FD)}]
\begin{enumerate}
\item
For all $t\in[0,T]$, $b(t,\cdot)$ is differentiable with continuous derivative and
\ce
\sup_{t\in[0,T],x\in\mR^d}|\nabla b(t,x)|<\infty;
\de
\item
for all $t\in[0,T]$, $\sigma(t,\cdot)$ is differentiable with continuous derivative and
\ce
\sup_{t\in[0,T],x\in\mR^d}|\nabla\sigma(t,x)|<\infty.
\de
\end{enumerate}
\item[\textbf{(FE)}]
$\sigma$ is bounded and uniformly nondegenerate, i.e., there exist two constants $C_{FB},C_{FE}>0$ such that $|\sigma(t,x)|\leq C_{FB}$ and $\sigma(t,x)\sigma^*(t,x)\geq C_{FE}I$.
\end{enumerate}

The following result can be found in \cite[Proposition 9]{BD2}.

\bt\label{FSDE existence}
\begin{enumerate}
\item
Under the assumption {\bf (FL)}, there exists a unique solution $\{X(\tau,t,x);\tau\in[t,T],x\in\mR^d\}$ to Eq. (\ref{SDE}).
\item
Moreover, for all $p\geq2$, there exists a constant $C_p>0$ only depending on the Lipschitz constants of $b$ and $\sigma$ such that for any $t\in[0,T]$,
\ce
\mE[\sup_{\tau\in[t,T]}|X(\tau,t,x)|^p]\leq C_p(1+|x|^p).
\de
\item
Furthermore, under the assumption {\bf (FD)}, the $L^2$-derivative of the solution to Eq. (\ref{SDE}) with respect to $x$, $\nabla X(\tau,t,x)$, exists and is the unique solution of the following matrix-valued SDE with jumps:
\ce
\nabla X(\tau,t,x)&=&I+\int_t^{\tau}\nabla b(s,X(s,t,x))\cdot \nabla X(s,t,x)ds\\
&&+\int_t^{\tau}\int_{\cO}\nabla\sigma(s,X(s-,t,x))u\cdot
\nabla X(s-,t,x)\widetilde{N}(ds,du).
\de
\item
Finally, the solution $\{X(\tau,t,x);\tau\in[t,T],x\in\mR^d\}$ belongs to $\cH_{\mD}^d$, and for every $t\in[0,T]$, the Malliavin derivative $X^{\sharp}(\tau,t,x)$ satisfies the following SDEs with jumps:
\ce
X^{\sharp}(\tau,t,x)&=&\int_t^{\tau}\nabla b(s,X(s,t,x))\cdot X^{\sharp}(s,t,x)ds\\
&&+\int_t^{\tau}\int_{\cO\times R}\sigma(s,X(s-,t,x))u^{\flat}N\odot\rho(ds,du,dr)\\
&&+\int_t^{\tau}\int_{\cO}\nabla\sigma(s,X(s-,t,x))u\cdot X^{\sharp}(s-,t,x)\widetilde{N}(ds,du).
\de
\end{enumerate}
\et

\subsection{Forward-Backward System}

We need the following assumptions on $\psi$ and $\phi$.

\begin{enumerate}
\item[\textbf{(BL)}]
The functions $\psi$ and $\phi$ are Lipschitz and with at most polynomial growth, i.e., there exist two constants $C_{BL}>0$ and $\mu\geq0$ such that, for all $t\in[0,T]$, $x\in\mR^d$, $y,y'\in\mR$, $z,z'\in\mR$,
\ce
|\psi(t,x,y,z)-\psi(t,x,y',z')|\leq C_{BL}(|y-y'|+|z-z'|),\\
|\psi(t,x,0,0)|\leq C_{BL}(1+|x|)^{\mu},\quad|\phi(x)|\leq C_{BL}(1+|x|)^{\mu}.
\de
\item[\textbf{(BD)}]
\begin{enumerate}
\item
$\psi(t,\cdot,\cdot,\cdot)\in\cC^{1,1,1}(\mR^d\times\mR\times\mR;\mR)$, $\phi\in\cC^1(\mR^d;\mR)$, for every $t\in[0,T]$;
\item
there exist two constants $C_{BD}>0$ and $m\geq0$ such that
\ce
|\nabla_x\psi(t,x,y,z)|\leq C_{BD}(1+|z|)(1+|x|+|y|)^m,\quad|\nabla_x\phi(x)|\leq C_{BD},
\de
for every $t\in[0,T]$, $x\in\mR^d$, $y\in\mR$ and $z\in\mR$.
\end{enumerate}
\end{enumerate}

The following results are important and we will give its proof in the Appendix A.

\bt\label{BSDE existence}
\begin{enumerate}
\item
Under the assumptions {\bf (FL)} and {\bf (BL)}, for any $p\geq1$, there exists a unique solution
\ce
(Y(\tau,t,x),Z(\tau,t,x,u))\in\cS^p\times\cM^p
\de
to Eq. (\ref{BSDE}).
\item
Moreover, for all $p\geq2$, there exists a constant $C_p>0$ only depending on the Lipschitz constants of $b$, $\sigma$, $\psi$ and $\phi$, such that, for any $t\in[0,T]$ and $x\in\mR^d$,
\ce
&&\mE[\sup_{\tau\in[t,T]}|Y(\tau,t,x)|^p+(\int_t^T\int_{\cO}|Z(\tau,t,x,u)|^2\nu(du)d\tau)^{\frac{p}{2}}]\\
&\leq&C_p\mE[\int_t^T|\psi(s,X(s,t,x),0,0)|^2ds]^{\frac{p}{2}}\\
&&+C_p\mE[|\phi(X_T)|^p],
\de
\ce
&&\mE[\sup_{\tau\in[t,T]}|Y(\tau,t,x)|^p+(\int_t^T\int_{\cO}|Z(\tau,t,x,u)|^2\nu(du)d\tau)^{\frac{p}{2}}]\\
&\leq&C_p(T-t)^{\frac{p}{2}}\mE[\int_t^T|\psi(s,X(s,t,x),Y(s,t,x),\int_{\cO}Z(s,t,x,u)l(u)\nu(du))|^2ds]^{\frac{p}{2}}\\
&&+C_p\mE[|\phi(X_T)|^p],
\de
and for any $\rho>0$,
\ce
&&\mE[\sup_{\tau\in[0,T]}e^{\rho\tau}|Y(\tau,t,x)|^2+\int_0^T\int_{\cO}e^{\rho\tau}|Z(\tau,t,x,u)|^2\nu(du)d\tau]\\
&\leq&\frac{C_2}{\rho}\mE[\int_0^Te^{\rho s}|\psi(s,X(s,t,x),Y(s,t,x),\int_{\cO}Z(s,t,x,u)l(u)\nu(du))|^2ds]\\
&&+C_2e^{\rho T}\mE[|\phi(X_T)|^2].
\de
\item
Furthermore, under the assumptions {\bf (FL)}, {\bf (BL)} and {\bf (BD)}, the $L^2$-derivative of the solution of Eq. (\ref{BSDE}) with respect to $x$,
\ce
(\nabla Y(\tau,t,x),\nabla Z(\tau,t,x,u)),
\de
exists and is the unique solution of the following backward SDEs with jumps:
\ce
&&\nabla Y(\tau,t,x)+\int_{\tau}^T\int_{\cO}\nabla Z(s-,t,x,u)\widetilde{N}(ds,du)\no\\
&=&-\int_{\tau}^T(\nabla_x\psi(s,X(s,t,x),Y(s,t,x),\int_{\cO}Z(s,t,x,u)l(u)\nu(du))\nabla X(s,t,x)\no\\
&&+\nabla_y\psi(s,X(s,t,x),Y(s,t,x),\int_{\cO}Z(s,t,x,u)l(u)\nu(du))\nabla Y(s,t,x)\no\\
&&+\nabla_z\psi(s,X(s,t,x),Y(s,t,x),\int_{\cO}Z(s,t,x,u)l(u)\nu(du))\int_{\cO}\nabla Z(s,t,x,u)l(u)\nu(du))ds\no\\
&&+\nabla\phi(X(T,t,x))\nabla X(T,t,x).
\de
\item
Finally, for every $t\in[0,T]$, the solution
\ce
\{Y(\tau,t,x),Z(\tau,t,x,u);\tau\in[t,T],x\in\mR^d,u\in\cO\}
\de
belongs to $\cH_{\mD}\times\cH_{\mD,\nu}$, and the Malliavin derivative
\ce
\{Y^{\sharp}(\tau,t,x),Z^{\sharp}(\tau,t,x,u);\tau\in[t,T],x\in\mR^d,u\in\cO\}
\de
satisfies the following backward SDEs with jumps:
\ce
&&Y^{\sharp}(\tau,t,x)+\int_{\tau}^T\int_{\cO}Z^{\sharp}(s-,t,x,u)\widetilde{N}(ds,du)\\
&=&-\int_{\tau}^T(\nabla_x\psi(s,X(s,t,x),Y(s,t,x),\int_{\cO}Z(s,t,x,u)l(u)\nu(du))X^{\sharp}(s,t,x)\\
&&+\nabla_y\psi(s,X(s,t,x),Y(s,t,x),\int_{\cO}Z(s,t,x,u)l(u)\nu(du))Y^{\sharp}(s,t,x)\\
&&+\nabla_z\psi(s,X(s,t,x),Y(s,t,x),\int_{\cO}Z(s,t,x,u)l(u)\nu(du))\int_{\cO}Z^{\sharp}(s,t,x,u)l(u)\nu(du))ds\\
&&+\nabla\phi(X(T,t,x))X^{\sharp}(T,t,x).
\de
\end{enumerate}
\et

\br
In fact, there has been relatively little research on the Malliavin differentiability of forward-backward SDEs with jumps before, and to our knowledge, only Delong in \cite{D} has considered this topic. When the Malliavin differentiability of the solution to Eq. (\ref{BSDE}) is obtained, we can even consider the problem of the existence and smoothness of the densities of forward-backward SDEs with jumps in our future work.
\er

\section{Bismut-Elworthy-Li Formulae}

We will make use of the following assumptions.

\begin{enumerate}
\item[\textbf{(L)}]
\begin{enumerate}
\item
The L\'{e}vy measure $\nu$ satisfies $\nu(\cO)=+\infty$.
\item
$k\in\cC_0^1(\cO;\mR^+\setminus\{0\})$, and there exists a positive constant $C_B$ such that
\ce
|\nabla\log k(u)|\leq C_B|u|^{-1},\quad u\in\cO.
\de
\item
there exists $\beta\in(0,2)$ such that
\begin{enumerate}
\item
For any $p\geq2$, there exists a positive constant $C_O$ such that for any $\varepsilon\in(0,1)$,
\ce
\int_{\{|u|\leq\varepsilon\}}|u|^p\nu(du)\leq C_O\varepsilon^{p-\beta};
\de
\item
For any $p\geq1$, there exists a positive constant $C_I$ such that for any $\tau\in[t,T]$ and $\varepsilon\in(0,1)$,
\ce
\mE[(\int_t^{\tau}\int_{\{|u|\leq\varepsilon\}}|u|^3N(ds,du))^{-p}]\leq C_I(((\tau-t)\varepsilon^{3-\beta})^{-p}+(\tau-t)^{-\frac{3p}{\beta}}).
\de
\end{enumerate}
\end{enumerate}
\end{enumerate}

Let us give some comments on theses assumptions. 
In \cite{T,WXZ1}, the following growth condition near the origin on the intensity measure is used:
\be\label{growth condition}
\lim_{\varepsilon\downarrow0}\varepsilon^{\beta-2}\int_{\{|u|<\varepsilon\}}|u|^2\nu(du)>0, ~\beta\in(0,2).
\ee
It can be seen from \cite[Lemma 2.5, Lemma 2.6]{WXZ1} that the assumption (\textbf{L}-c) is satisfied under this growth condition. It is easy to see that if $\nu(du)/du=k(u)=a(u)|u|^{-d-\beta}$ with
\ce
a(u)=a(-u),\quad 0<a_0\leq a(u)\leq a_1,\quad|\nabla a(u)|\leq a_2,
\de
then the growth condition (\ref{growth condition}) and the assumption (\textbf{L}-a,b) are satisfied. The L\'{e}vy processes corresponding to this $\nu$ is called $\beta$-stable-like process as in \cite{CK}.

However, Proposition 4.11 in \cite{RZ} reveals that the growth condition (\ref{growth condition}) is stronger than the following one: there exists $\beta\in(0,2)$ such that
\ce
(0,\infty)\ni\kappa=\limsup_{\lambda\rightarrow\infty}\frac{1}{\lambda^{\beta/3}}\int_{\cO}(1-\exp\{-\lambda|u|^3\})\nu(du).
\de
The following proposition, whose detailed proof  is included in the Appendix B,  shows that the assumption (\textbf{L}-c) still holds under this relatively weak condition.
 
\bp\label{weak condition}
If there exists $\beta\in(0,2)$ such that
\be\label{inverse}
(0,\infty)\ni\kappa=\limsup_{\lambda\rightarrow\infty}\frac{1}{\lambda^{\beta/3}}\int_{\cO}(1-\exp\{-\lambda|u|^3\})\nu(du),
\ee
then
\begin{enumerate}
\item
for any $p\geq2$, there exists a positive constant $C_O$ such that for any $\varepsilon\in(0,1)$,
\ce
\int_{\{|u|\leq\varepsilon\}}|u|^p\nu(du)\leq C_O\varepsilon^{p-\beta};
\de
\item
for any $p\geq1$, there exists a positive constant $C_I$ such that for any $\tau\in[t,T]$ and $\varepsilon\in(0,1)$,
\ce
\mE[(\int_t^{\tau}\int_{\{|u|\leq\varepsilon\}}|u|^3N(ds,du))^{-p}]\leq C_I(((\tau-t)\varepsilon^{3-\beta})^{-p}+(\tau-t)^{-\frac{3p}{\beta}}).
\de
\end{enumerate}
\ep

Fix $\varepsilon\in(0,1)$ and let $\zeta_{\varepsilon}(u)$ be a smooth real function with
\ce
\zeta_{\varepsilon}(u)=
\begin{cases}
|u|^3,&|u|\leq\frac{\varepsilon}{3},\\
0,&|u|>\frac{2\varepsilon}{3}
\end{cases}
\de
and
\ce
|\nabla\zeta_{\varepsilon}(u)|\leq C|u|^2,
\de
where $C$ is independence of $\varepsilon$. Denote by $\cC_0^\infty(\cO)$ the set of $\cC^\infty$-functions defined on $\cO$ and with compact support and by $H$ the subspace consisting of $f\in L^2(\nu)\cap L^1(\nu)$ such that $f\in\cC_0^{\infty}(\cO)$. Then the bilinear form
\ce
e_{\varepsilon}(\phi,\psi)=\frac{1}{2}\sum_{i=1}^l\int_{\cO}\zeta_{\varepsilon}(u)\partial_i\phi(u)\partial_i\psi(u)k(u)du,\quad\forall\phi,\psi\in H
\de
is closable and its closure, which will be  denoted by $(\mathbf{d}_{\varepsilon},e_{\varepsilon})$, is a local symmetric Dirichlet form on $L^2(kdu)$ which admits a carr\'{e} du champ operator $\gamma_{\varepsilon}$ given by
\ce
\quad\gamma_{\varepsilon}[\phi,\psi]=\sum_{i=1}^l\zeta_{\varepsilon}\partial_i\phi\partial_i\psi,\quad\forall \phi,\psi\in\mathbf{d}_{\varepsilon}.
\de

Since
\ce
\gamma_{\varepsilon}[\phi](u)=\zeta_{\varepsilon}(u)\sum_{i=1}^l(\partial_i\phi(u))^2,\quad\forall \phi\in\cC_0^{\infty}(\cO),
\de
and the identity map $j_{\varepsilon}$ belongs to $\mathbf{d}_{\varepsilon}$, we have $\gamma_{\varepsilon}[j_{\varepsilon},j_{\varepsilon}^*](u)=\zeta_{\varepsilon}(u)I$. Therefore,
\ce
j_{\varepsilon}^{\flat}(u,r)=\zeta_{\varepsilon}^{1/2}(u)\xi(r).
\de
Here
\ce
\xi=
\left(\begin{array}{c}
\xi_1(r) \\
\xi_2(r) \\
\vdots \\
\xi_l(r)
\end{array}\right),
\de
where $\xi_i\in L^2(\rho)$, $i=1,2,\cdots,l$, satisfy
\ce
\int_R\xi_i(r)\rho(dr)=0
\de
and
\ce
\int_R\xi_i(r)\xi_j(r)\rho(dr)=\delta_{ij}
\de
for all $i,j=1,\cdots,l$. Consequently,
\ce
\phi^{\flat}(u,r)=\zeta_{\varepsilon}^{1/2}(u)\sum_{i=1}^l\partial_i\phi(u)\xi_i(r).
\de
Here and in the sequel, the operators $\flat$, $\delta_\flat$, $\sharp$ and $\delta_\sharp$ corresponding to $(d_\varepsilon, e_\varepsilon)$ depend of course on $\varepsilon$, but for notational simplicity, we will omit this $\varepsilon$ as a subscript. For a $U(u,r)\in L^2(\nu\times \rho)$ which is $\cC^1$ in the first variable,
\be\label{theta}
[\delta_\flat U](u)=-\int_R[\theta_{\varepsilon}(u,r)U(u,r)+\zeta_{\varepsilon}^{1/2}(u)\sum_{i=1}^l\partial_iU(u,r)\xi_i(r)]\rho(dr),
\ee
where
\be
\begin{split}
\theta_{\varepsilon}(u,r)&=\sum_{i=1}^l\frac{\partial_i(\zeta_{\varepsilon}^{1/2}(u)k(u))}{k(u)}\xi_i(r)\\
&=\sum_{i=1}^l(\partial_i\log k(u)\zeta_{\varepsilon}^{1/2}(u)+\frac12\zeta_{\varepsilon}^{-1/2}(u)\partial_i\zeta_{\varepsilon}(u))\xi_i(r).
\end{split}
\ee
In fact, according to the integration by parts formula and Fubini theorem, we have
\ce
&&\int_{\cO}\int_RU(u,r)\phi^{\flat}(u,r)\rho(dr)\nu(du)\\
&=&\int_R\int_{\cO}U(u,r)\zeta_{\varepsilon}^{1/2}(u)\sum_{i=1}^l\partial_i\phi(u)\xi_i(r)k(u)du\rho(dr)\\
&=&\int_R\sum_{i=1}^l\xi_i(r)(\int_{\cO}\partial_i\phi(u)U(u,r)\zeta_{\varepsilon}^{1/2}(u)k(u)du)\rho(dr)\\
&=&\int_R\sum_{i=1}^l\xi_i(r)(\int_{\cO}\phi(u)\partial_i(U(u,r)\zeta_{\varepsilon}^{1/2}(u)k(u))du)\rho(dr)\\
&=&\int_R\sum_{i=1}^l\xi_i(r)(\int_{\cO}\phi(u)\partial_i(\zeta_{\varepsilon}^{1/2}(u)k(u))U(u,r)du)\rho(dr)\\
&&+\int_R\sum_{i=1}^l\xi_i(r)(\int_{\cO}\phi(u)\zeta_{\varepsilon}^{1/2}(u)k(u)\partial_iU(u,r)du)\rho(dr)\\
&=&\int_{\cO}(\int_R(\sum_{i=1}^l\frac{\partial_i(\zeta_{\varepsilon}^{1/2}(u)k(u))}{k(u)}\xi_i(r)\\
&&+\zeta_{\varepsilon}^{1/2}(u)\sum_{i=1}^l\partial_iU(u,r)\xi_i(r))\rho(dr))\phi(u)\nu(du),
\de
then the desired result can be obtained through the definition of $[\delta_\flat U](u)$.

Now the derivative formulae of Bismut-Elworthy-Li's type for forward-backward SDEs with respect to Poisson random measure can be established as follows.

\bt\label{Bismut formula}
Suppose that the assumptions {\bf (L)}, {\bf (FD)}, {\bf (FE)}, {\bf (BL)} and {\bf (BD)} are satisfied. Then for any $\tau\in[t,T]$ and $x,h\in\mR^d$, we have
\ce
&&\mE[\nabla Y(\tau,t,x)h]\\
&=&-\mE[\int_{\tau}^T\psi(s,X(s,t,x),Y(s,t,x),\int_{\cO}Z(s,t,x,u)l(u)\nu(du))\cdot[\delta_{\sharp}(F_{\varepsilon}^*(s)\cdot h)\cdot G_{\varepsilon}^{-1}(s)\\
&&-\widehat{\mE}[F_{\varepsilon}^*(s)\cdot h\cdot (G_{\varepsilon}^{-1}(s))^{\sharp}]]ds]\\
&&+\mE[\phi(X(T,t,x))\cdot[\delta_{\sharp}(F_{\varepsilon}^*(T)\cdot h)\cdot G_{\varepsilon}^{-1}(T)\\
&&-\widehat{\mE}[F_{\varepsilon}^*(T)\cdot h\cdot (G_{\varepsilon}^{-1}(T))^{\sharp}]]],
\de
where
\ce
F_{\varepsilon}(\tau)&=&\int_t^{\tau}\int_{\cO}\int_R\sigma^{-1}(\alpha,X(\alpha-,t,x))\cdot(I+\nabla\sigma(\alpha,X(\alpha-,t,x))y)\\
&&\cdot\nabla X(\alpha-,t,x)\cdot j_{\varepsilon}^{\flat}(y,r)N\odot\rho(d\alpha,dy,dr)
\de
and
\ce
G_{\varepsilon}(\tau)=\int_t^{\tau}\int_{\cO}\zeta_{\varepsilon}(y)N(d\alpha,dy).
\de

\et

The following gradient estimate will play an important role in the subsequent section.

\bt\label{gradient estimate}
Suppose that $\beta\in(1,2)$ and the assumptions {\bf (L)}, {\bf (FD)}, {\bf (FE)}, {\bf (BL)} and {\bf (BD)} are satisfied. Then there exists a positive constant $C$, which is independent of the constants $C_{BL}$ and $C_{BD}$, such that for all $t\in[0,T]$ and $x\in\mR^d$,
\ce
|\nabla Y(t,t,x)|\leq C(T-t)^{-1/\beta}(1+|x|)^{\mu}.
\de
\et

\br
According to Theorem \ref{Bismut formula}, it is worth mentioning that the estimate on $\nabla Y(t,t,x)$ does not depend on the constant in the assumption {\bf (BD)}, and hence the regularity assumption of $\psi$ and $\phi$ required in the assumption {\bf (BD)} can be removed, and this is the starting point for our application in Section 5.
\er

\br
Li in \cite{L} has considered the mean-field forward and backward SDEs with jumps and associated nonlocal quasi-linear integral-PDEs of mean-field type. In the future work, using the lent particle method, we will attempt to obtain the Bismut-Elworthy-Li formule for the mean-field forward-backward SDEs with jumps and apply it to the quasi-linear integral-PDEs of mean field type to weaken the conditions of initial datum and coefficients.
\er

\br
Although our focus has been exclusively on the jump case in this paper, we note that the results presented here can be extended to incorporate Brownian motion as well. However, to maintain the conciseness and focus of the current work, we have omitted this extension. The necessary modifications to include Brownian motion are relatively straightforward and do not present significant additional difficulties in the proofs.
\er

Now we are in a position to prove Theorem \ref{Bismut formula} and Theorem \ref{gradient estimate}.

\subsection{Some Useful Lemmas}

To prove Theorem \ref{Bismut formula} and Theorem \ref{gradient estimate}, we need several lemmas. The first one is the following Kunita's inequalities which can be found in \cite[Theorem 2.11]{K}.

\bl\label{bdg}
\begin{enumerate}
\item
For any $p\geq1$, there is a positive constant $C$ such that
\be\label{bdg1}
\begin{split}
&\mE[\sup_{t\in[0,T]}|\int_0^t\int_{\cO}\psi(s,u)N(ds,du)|^p]\\
\leq&C\mE[\int_0^T\int_{\cO}|\psi(s,u)|\nu(du)ds]^p+C\mE[\int_0^T\int_{\cO}|\psi(s,u)|^p\nu(du)ds],
\end{split}
\ee
where $\psi(s,u):\Omega\times[0,T]\times\mR^l\rightarrow\mR^d$ is a predictable process such that the term on the right-hand side of (\ref{bdg1}) is finite.
\item
For any $p\geq2$, there is a positive constant $C$ such that
\be\label{bdg2}
\begin{split}
&\mE[\sup_{t\in[0,T]}|\int_0^t\int_{\cO}\psi(s,u)\widetilde{N}(ds,du)|^p]\\
\leq&C\mE[\int_0^T\int_{\cO}|\psi(s,u)|^2\nu(du)ds]^{p/2}+C\mE[\int_0^T\int_{\cO}|\psi(s,u)|^p\nu(du)ds],
\end{split}
\ee
where $\psi(s,u):\Omega\times[0,T]\times\mR^l\rightarrow\mR^d$ is a predictable process such that the term on the right-hand side of (\ref{bdg2}) is finite.
\end{enumerate}
\el

\bl\label{flow}
Under the assumptions {\bf (FD)}, {\bf (BL)} and {\bf (BD)}, $\mP$-a.s., for a.e. $s\in[t,T]$ and $u\in\cO$, we have
\ce
&&(\nabla Y(\tau,s,X(s,t,x))\nabla X(s,t,x),\nabla Z(\tau,s,X(s,t,x),u)\nabla X(s,t,x))\\
&=&(\nabla Y(\tau,t,x),\nabla Z(\tau,t,x,u)).
\de
\el

\begin{proof}
First of all, using the uniqueness of the solution of Eq. (\ref{SDE}), for any $t\leq s\leq\tau\leq T$, we have
\ce
\mP-\text{a.s.},\quad X(\tau,s,X(s,t,x))=X(\tau,t,x).
\de
It is also easy to see that for any $h\in\mR^d$ and $t\leq s\leq\tau\leq T$, we have
\ce
\mP-\text{a.s.},\quad\nabla X(\tau,s,X(s,t,x))\nabla X(s,t,x)h=\nabla X(\tau,t,x)h.
\de
Since the solution of the backward SDE with jumps is uniquely determined, on an interval $[s,T]$, by the values of the process $X$ on the same interval, for $t\leq s\leq\tau\leq T$, we have $\mP$-a.s.,
\ce
Y(\tau,s,X(s,t,x))=Y(\tau,t,x),\\
Z(\tau,s,X(s,t,x),u)=Z(\tau,t,x,u),
\de
and hence the desired result can be obtained by the chain rule.
\end{proof}

\bl\label{derivative Malliavin}
Under the assumptions {\bf (FD)}, {\bf (BL)} and {\bf (BD)}, for all $\tau\in[t,T]$, $x\in\mR^d$ and $u\in\cO$, we have
\ce
X^{\sharp}(\tau,t,x)&=&\int_t^{\tau}\int_{\cO}\int_R\nabla X(\tau,\alpha-,x)(\sigma(\alpha,X(\alpha-,t,x))y)^{\flat}N\odot\rho(d\alpha,dy,dr),\\
Y^{\sharp}(\tau,t,x)&=&\int_t^{\tau}\int_{\cO}\int_R\nabla Y(\tau,\alpha-,x)(\sigma(\alpha,X(\alpha-,t,x))y)^{\flat}N\odot\rho(d\alpha,dy,dr),\\
Z^{\sharp}(\tau,t,x,u)&=&\int_t^{\tau}\int_{\cO}\int_R\nabla Z(\tau,\alpha-,x,u)(\sigma(\alpha,X(\alpha-,t,x))y)^{\flat}N\odot\rho(d\alpha,dy,dr).
\de
\el

\begin{proof}
For $(\alpha,y,u)\in[t,T]\times\mR^d\times\cO$, setting
\ce
X^{(\alpha,y)}(\tau,t,x):=\varepsilon^+_{(\alpha,y)}X(\tau,t,x),\quad\nabla X^{\alpha,y}(\tau,\alpha,x):=\varepsilon^+_{(\alpha,y)}\nabla X(\tau,\alpha,x),\\
Y^{(\alpha,y)}(\tau,t,x):=\varepsilon^+_{(\alpha,y)}Y(\tau,t,x),\quad\nabla Y^{\alpha,y}(\tau,\alpha,x):=\varepsilon^+_{(\alpha,y)}\nabla Y(\tau,\alpha,x),\\
Z^{(\alpha,y)}(\tau,t,x,u):=\varepsilon^+_{(\alpha,y)}Z(\tau,t,x,u),\quad\nabla Z^{\alpha,y}(\tau,\alpha,x,u):=\varepsilon^+_{(\alpha,y)}\nabla Z(\tau,\alpha,x,u).
\de
It is obvious that
\ce
&&X^{(\alpha,y)}(\tau,t,x)\\
&=&x+\int_t^{\alpha}b(s,X(s,t,x))ds+\int_t^{\alpha}\int_{\cO}\sigma(s,X(s-,t,x))u\widetilde{N}(ds,du)\\
&&+\int_{\alpha}^{\tau}b(s,X^{(\alpha,y)}(s,t,x))ds\\
&&+\sigma(\alpha,X(\alpha-,t,x))y+\int_{\alpha}^{\tau}\int_{\cO}\sigma(s,X^{(\alpha,y)}(s-,t,x))u\widetilde{N}(ds,du).
\de
Using Theorem \ref{FSDE existence}, we also have
\ce
&&\nabla X^{(\alpha,y)}(\tau,\alpha,x)\\
&=&I+\int_{\alpha}^{\tau}\nabla b(s,X^{(\alpha,y)}(s,t,x))\cdot\nabla X^{(\alpha,y)}(s,\alpha,x)ds\\
&&+\int_{\alpha}^{\tau}\int_{\cO}\nabla \sigma(s,X^{(\alpha,y)}(s-,t,x))u\cdot\nabla X^{(\alpha,y)}(s-,\alpha,x)\widetilde{N}(ds,du).
\de
Taking the gradient of $X^{(\alpha,y)}(t,x)$ with respect to the variable $y$, we obtain
\ce
&&(X^{(\alpha,y)}(\tau,t,x))^{\flat}\\
&=&\int_{\alpha}^{\tau}\nabla b(s,X^{(\alpha,y)}(s,t,x))\cdot(X^{(\alpha,y)}(s-,t,x))^{\flat}ds\\
&&+(\sigma(\alpha,X(\alpha-,t,x))y)^{\flat}+\int_{\alpha}^{\tau}\int_{\cO}\nabla\sigma(s,X^{(\alpha,y)}(s-,t,x))u\cdot(X^{(\alpha,y)}(s-,t,x))^{\flat}\widetilde{N}(ds,du).
\de
Then for $(\alpha,y)\in[0,T]\times\mR^d$, we have the following relation:
\be\label{relation}
(X^{(\alpha,y)}(\tau,t,x))^{\flat}&=&\nabla X^{(\alpha,y)}(\tau,\alpha,x)(\sigma(\alpha,X(\alpha-,t,x))y)^{\flat}.
\ee
Now, we calculate gradient on the bottom Dirichlet structure and then we take back the particle:
\ce
\forall\tau\in[t,T],\quad\varepsilon_{(\alpha,y)}^-(X^{(\alpha,y)}(\tau,t,x))^{\flat}=\nabla X(\tau,\alpha,x)(\sigma(\alpha,X(\alpha-,t,x))y)^{\flat},
\de
which is an equality $\mP_{N\odot\rho}$-a.e. Finally by the formula of gradient operator (\ref{gradient}), we have
\ce
X^{\sharp}(\tau,t,x)=\int_t^{\tau}\int_{\cO}\int_R\nabla X(\tau,\alpha,x)(\sigma(\alpha,X(\alpha-,t,x))y)^{\flat}N\odot\rho(d\alpha,dy,dr).
\de

On the other hand, using the same argument as above, we obtain
\ce
&&Y^{(\alpha,y)}(\tau,t,x)+\int_{\tau}^T\int_{\cO}Z^{(\alpha,y)}(s-,t,x,u)\widetilde{N}(ds,du)+Z(\alpha,t,x,y)1_{\{\tau\in[0,\alpha]\}}\\
&=&-\int_{\tau}^T\psi(s,X^{(\alpha,y)}(s,t,x),Y^{(\alpha,y)}(s,t,x),\int_{\cO}Z^{(\alpha,y)}(s,t,x,u)l(u)\nu(du))ds\\
&&+\nabla\phi(X^{(\alpha,y)}(T,t,x)).
\de
Then for any $\tau\in[\alpha,T]$, we have
\be\label{derivative bottom}
\begin{split}
&(Y^{(\alpha,y)}(\tau,t,x))^{\flat}+\int_{\tau}^T\int_{\cO}(Z^{(\alpha,y)}(s-,t,x,u))^{\flat}\widetilde{N}(ds,du)\\
=&-\int_{\tau}^T(\nabla_x\psi(s,X^{(\alpha,y)}(s,t,x),Y^{(\alpha,y)}(s,t,x),\\
&\int_{\cO}Z^{(\alpha,y)}(s,t,x,u)l(u)\nu(du))(X^{(\alpha,y)}(s,t,x))^{\flat}\\
&+\nabla_y\psi(s,X^{(\alpha,y)}(s,t,x),Y^{(\alpha,y)}(s,t,x),\\
&\int_{\cO}Z^{(\alpha,y)}(s,t,x,u)l(u)\nu(du))(Y^{(\alpha,y)}(s,t,x))^{\flat}\\
&+\nabla_z\psi(s,X^{(\alpha,y)}(s,t,x),Y^{(\alpha,y)}(s,t,x),\\
&\int_{\cO}Z^{(\alpha,y)}(s,t,x,u)l(u)\nu(du))\int_{\cO}(Z^{(\alpha,y)}(s,t,x,u))^{\flat}l(u)\nu(du))ds\\
&+\nabla\phi(X^{(\alpha,y)}(T,t,x))(X^{(\alpha,y)}(T,t,x))^{\flat}.
\end{split}
\ee
Applying the similar procedure as in the proof of Theorem \ref{BSDE existence}, for any $\alpha\in[0,\tau]$ and $y,h\in\mR^d$, it holds that
\ce
&&\nabla Y^{(\alpha,y)}(\tau,\alpha,x)h+\int_{\tau}^T\int_{\cO}\nabla Z^{(\alpha,y)}(s-,\alpha,x,u)h\widetilde{N}(ds,du)\no\\
&=&-\int_{\tau}^T(\nabla_x\psi(s,X^{(\alpha,y)}(s,\alpha,x),Y^{(\alpha,y)}(s,\alpha,x),\int_{\cO}Z^{(\alpha,y)}(s,\alpha,x,u)l(u)\nu(du))\nabla X^{(\alpha,y)}(s,\alpha,x)h\no\\
&&+\nabla_y\psi(s,X^{(\alpha,y)}(s,\alpha,x),Y^{(\alpha,y)}(s,\alpha,x),\int_{\cO}Z^{(\alpha,y)}(s,\alpha,x,u)l(u)\nu(du))\nabla Y^{(\alpha,y)}(s,\alpha,x)h\no\\
&&+\nabla_z\psi(s,X^{(\alpha,y)}(s,\alpha,x),Y^{(\alpha,y)}(s,\alpha,x),\int_{\cO}Z^{(\alpha,y)}(s,\alpha,x,u)l(u)\nu(du))\\
&&\int_{\cO}\nabla Z^{(\alpha,y)}(s,t,x,u)l(u)\nu(du)h)ds\no\\
&&+\nabla\phi(X^{(\alpha,y)}(T,\alpha,x))\nabla X^{(\alpha,y)}(T,\alpha,x)h.
\de
Now given $v\in\mR^d$ and $t\in[0,\alpha]$, we can replace $x$ by $X(\alpha,t,x)$ and $h$ by $\sigma(\alpha,X(\alpha-,t,x))y^{\flat}v$ in this equation, and hence using (\ref{relation}) we have
\ce
&&\nabla Y^{(\alpha,y)}(\tau,\alpha,X(\alpha,t,x))\sigma(\alpha,X(\alpha,t,x))y^{\flat}v\\
&&+\int_{\tau}^T\int_{\cO}\nabla Z^{(\alpha,y)}(s-,\alpha,X(\alpha,t,x),u)\sigma(\alpha,X(\alpha,t,x))y^{\flat}v\widetilde{N}(ds,du)\no\\
&=&-\int_{\tau}^T(\nabla_x\psi(s,X^{(\alpha,y)}(s,\alpha,X(\alpha,t,x)),Y^{(\alpha,y)}(s,\alpha,X(\alpha,t,x)),\no\\
&&\int_{\cO}Z^{(\alpha,y)}(s,\alpha,X(\alpha,t,x),u)l(u)\nu(du))(X^{(\alpha,y)}(s,t,x))^{\flat}v\no\\
&&+\nabla_y\psi(s,X^{(\alpha,y)}(s,\alpha,X(\alpha,t,x)),Y^{(\alpha,y)}(s,\alpha,X(\alpha,t,x)),\no\\
&&\int_{\cO}Z^{(\alpha,y)}(s,\alpha,X(\alpha,t,x),u)l(u)\nu(du))\nabla Y^{(\alpha,y)}(s,\alpha,X(\alpha,t,x))\sigma(\alpha,X(\alpha-,t,x))y^{\flat}v\no\\
&&+\nabla_z\psi(s,X^{(\alpha,y)}(s,\alpha,X(\alpha,t,x)),Y^{(\alpha,y)}(s,\alpha,X(\alpha,t,x)),\no\\
&&\int_{\cO}Z^{(\alpha,y)}(s,\alpha,X(\alpha,t,x),u)l(u)\nu(du))\no\\
&&\int_{\cO}\nabla Z^{(\alpha,y)}(s,t,X(\alpha,t,x),u)l(u)\nu(du)\sigma(\alpha,X(\alpha-,t,x))y^{\flat}v)ds\no\\
&&+\nabla\phi(X^{(\alpha,y)}(T,t,x))(X^{(\alpha,y)}(T,t,x))^{\flat}v.
\de
This implies by (\ref{derivative bottom}) and Lemma \ref{flow} that
\ce
(Y^{(\alpha,y)}(\tau,t,x))^{\flat}&=&\nabla Y^{(\alpha,y)}(\tau,\alpha,X(\alpha,t,x))(\sigma(\alpha,X(\alpha-,t,x))y)^{\flat}\\
&=&\nabla Y^{(\alpha,y)}(\tau,\alpha,x)(\sigma(\alpha,X(\alpha-,t,x))y)^{\flat}
\de
and
\ce
(Z^{(\alpha,y)}(\tau,t,x,u))^{\flat}&=&\nabla Z^{(\alpha,y)}(\tau,\alpha,X(\alpha,t,x),u)(\sigma(\alpha,X(\alpha-,t,x))y)^{\flat}\\
&=&\nabla Z^{(\alpha,y)}(\tau,\alpha,x,u)(\sigma(\alpha,X(\alpha-,t,x))y)^{\flat}.
\de
Then the desired result can be obtained by the formula of gradient operator (\ref{gradient}).
\end{proof}

\bl\label{derivative initial}
Under the assumptions {\bf (L)}, {\bf (FD)}, {\bf (FE)}, {\bf (BL)} and {\bf (BD)} we have
\ce
\nabla_xX(\tau,t,x)&=&\widehat{\mE}[(X(\tau,t,x))^{\sharp}F_{\varepsilon}^*(\tau)G_{\varepsilon}^{-1}(\tau)],\\
\nabla_xY(\tau,t,x)&=&\widehat{\mE}[(Y(\tau,t,x))^{\sharp}F_{\varepsilon}^*(\tau)G_{\varepsilon}^{-1}(\tau)],\\
\nabla_xZ(\tau,t,x,u)&=&\widehat{\mE}[(Z(\tau,t,x,u))^{\sharp}F_{\varepsilon}^*(\tau)G_{\varepsilon}^{-1}(\tau)],
\de
where $F_{\varepsilon}$ and $G_{\varepsilon}$ are defined as in Theorem \ref{Bismut formula}.
\el

\begin{proof}
First of all, by Lemma \ref{derivative Malliavin}, we have
\ce
X^{\sharp}(\tau,t,x)=\int_t^{\tau}\int_{\cO}\int_R\nabla X(\tau,\alpha,x)\cdot(\sigma(\alpha,X(\alpha-,t,x))u)^{\flat}N\odot\rho(d\alpha,du,dr).
\de
On the other hand we have
\ce
\nabla_xX(\alpha,t,x)-\nabla_xX(\alpha-,t,x)=\nabla \sigma(\alpha,X(\alpha-,t,x))\cdot\Delta P(\alpha)\cdot \nabla_xX(\alpha-,t,x),
\de
where
\ce
P(s)=\int_t^s\int_{\cO}u\widetilde{N}(ds,du),
\de
or equivalently,
\ce
(I+\nabla\sigma(\alpha,X(\alpha-,t,x))y)\cdot \nabla_xX(\alpha-,t,x)=\nabla_xX(\alpha,t,x).
\de
Then by Proposition \ref{marked Poisson measure} and Lemma \ref{flow}, we have
\ce
&&\widehat{\mE}[X^{\sharp}(\tau,t,x)\cdot F_{\varepsilon}^*(\tau)\cdot h\cdot G_{\varepsilon}^{-1}(\tau)]\\
&=&\widehat{\mE}[\int_t^{\tau}\int_{\cO}\int_R\nabla_xX(\tau,\alpha,x)\cdot(\sigma(\alpha,X(\alpha-,t,x))u)^{\flat}N\odot\rho(d\alpha,du,dr)\\
&&F_{\varepsilon}^*(\tau)\cdot h\cdot G_{\varepsilon}^{-1}(\tau)]\\
&=&\widehat{\mE}[\int_t^{\tau}\int_{\cO}\int_R\nabla_xX(\tau,\alpha,x)\cdot\sigma(\alpha,X(\alpha-,t,x))\cdot j_{\varepsilon}^{\flat}(u,r)N\odot\rho(d\alpha,du,dr)\\
&& F_{\varepsilon}^*(\tau)\cdot h\cdot G_{\varepsilon}^{-1}(\tau)]\\
&=&\int_t^{\tau}\int_{\cO}\nabla_xX(\tau,\alpha,x)\cdot\sigma(\alpha,X(\alpha-,t,x))\cdot\gamma_{\varepsilon}[j,j^*](u)\\
&&\cdot\sigma^{-1}(\alpha,X(\alpha-,t,x))\cdot(I+\nabla\sigma(\alpha,X(\alpha-,t,x))u)\cdot \nabla_xX(\alpha-,t,x)N(d\alpha,du)\cdot h\cdot G_{\varepsilon}^{-1}(\tau)\\
&=&\nabla_xX(\tau,t,x)\cdot h\cdot\int_t^{\tau}\int_{\cO}\zeta_{\varepsilon}(u)N(d\alpha,du)\cdot G_{\varepsilon}^{-1}(\tau)\\
&=&\nabla_xX(\tau,t,x)\cdot h.
\de
Using Lemma \ref{flow} and Lemma \ref{derivative Malliavin}, the representation of $\nabla_xY(\tau,t,x)$ and $\nabla_xZ(\tau,t,x,u)$ can be obtained similarly. Then we complete the proof.
\end{proof}

\bl\label{divergence operator}
Let $F_{\varepsilon}(\tau)$, $G_{\varepsilon}(\tau)$ be defined as in Theorem \ref{Bismut formula}. Then for any $\tau\in[t,T]$ and $h\in\mR^d$, we have
\ce
&&\delta_{\sharp}(F_{\varepsilon}^*(\tau)\cdot h)\\
&=&\int_t^{\tau}\int_{\cO}\int_{R}\widehat{\mE}[-\theta_{\varepsilon}(y,r)\cdot(\sigma^{-1}(\alpha,X(\alpha-,t,x))\cdot(I+\nabla\sigma(\alpha,X(\alpha-,t,x))y)\\
&&\cdot\nabla X(\alpha-,t,x)\cdot j_{\varepsilon}^{\flat}(y,r))^*\cdot h\\
&&-\zeta_{\varepsilon}^{1/2}(y)\cdot\sum_{i=1}^l(\sigma^{-1}(\alpha,X(\alpha-,t,x))\cdot\nabla\sigma(\alpha,X(\alpha-,t,x))e_i\cdot\nabla X(\alpha-,t,x)\cdot j_{\varepsilon}^{\flat}(y,r)\\
&&+\sigma^{-1}(\alpha,X(\alpha-,t,x))\cdot(I+\nabla\sigma(\alpha,X(\alpha-,t,x))y)\cdot\nabla X(\alpha-,t,x)\cdot\partial_i\zeta_{\varepsilon}^{1/2}(y)\cdot\xi(r))^*\cdot h]\\
&&\rho(dr)N(d\alpha,dy),
\de
where $e_i=(0,\cdots,0,1,0,\cdots,0)^*$, and
\ce
G_{\varepsilon}^{\sharp}(\tau)=\int_t^{\tau}\int_{\cO}\int_R\zeta_{\varepsilon}^{1/2}(y)\sum_{i=1}^l\partial_i\zeta_{\varepsilon}(y)\xi_i(r)N\odot\rho(d\alpha,dy,dr).
\de
\el

\begin{proof}
First of all, the second conclusion can be easily obtained by \cite[Proposition 8]{BD2}.

Now we prove the first one. Let $(\alpha,y,r)\in[0,T]\times\mR^l\times R$. Then it is obvious that
\ce
&&\varepsilon_{(\alpha,y,r)}^+(F_{\varepsilon}^*(\tau)\cdot h)\\
&=&(\int_t^{\tau}\int_{\cO}\int_R\sigma^{-1}(\alpha,\varepsilon_{(\alpha,y)}^+X(\alpha-,t,x))\cdot(I+\nabla\sigma(\alpha,\varepsilon_{(\alpha,y)}^+X(\alpha-,t,x))y)\\
&&\cdot\varepsilon_{(\alpha,y)}^+(\nabla X(\alpha-,t,x))\cdot j_{\varepsilon}^{\flat}(y,r)N\odot\rho(d\alpha,dy,dr))^*\cdot h\\
&&+(\sigma^{-1}(\alpha,X(\alpha-,t,x))\cdot(I+\nabla\sigma(\alpha,X(\alpha-,t,x))y)\cdot\nabla X(\alpha-,t,x)\cdot j_{\varepsilon}^{\flat}(y,r))^*\cdot h.
\de
Therefore we have
\ce
&&\partial_i(\varepsilon_{(\alpha,y,r)}^+(F_{\varepsilon}^*(\tau)\cdot h))\\
&=&((\sigma^{-1}(\alpha,X(\alpha-,t,x))\cdot\nabla\sigma(\alpha,X(\alpha-,t,x))e_i\cdot \nabla X(\alpha-,t,x)\cdot j_{\varepsilon}^{\flat}(y,r)\\
&&+\sigma^{-1}(\alpha,X(\alpha-,t,x))\cdot(I+\nabla\sigma(\alpha,X(\alpha-,t,x))y)\cdot \nabla X(\alpha-,t,x)\cdot \partial_ij_{\varepsilon}^{\flat}(y,r))^*\cdot h\\
&=&(\sigma^{-1}(\alpha,X(\alpha-,t,x))\cdot\nabla\sigma(\alpha,X(\alpha-,t,x))e_i\cdot \nabla X(\alpha-,t,x)\cdot j_{\varepsilon}^{\flat}(y,r)\\
&&+\sigma^{-1}(\alpha,X(\alpha-,t,x))\cdot(I+\nabla\sigma(\alpha,X(\alpha-,t,x))y)\cdot \nabla X(\alpha-,t,x)\cdot\partial_i\zeta_{\varepsilon}^{1/2}(y)\cdot\xi(r))^*\cdot h.
\de
By (\ref{theta}) we obtain
\ce
&&\widehat{\mE}[\delta_{\flat}(\varepsilon_{(\alpha,y,r)}^+(F_{\varepsilon}^*(\tau)\cdot h))]\\
&=&\int_{R}\widehat{\mE}[-\theta_{\varepsilon}(y,r)\cdot\varepsilon_{(\alpha,y,r)}^+(F_{\varepsilon}^*(\tau)\cdot h)-\zeta_{\varepsilon}^{1/2}(y)\cdot\sum_{i=1}^l\partial_i(\varepsilon_{(\alpha,y,r)}^+(F_{\varepsilon}^*(\tau)\cdot h))\cdot\xi(r)]\rho(dr)\\
&=&\int_{R}\widehat{\mE}[-\theta_{\varepsilon}(y,r)\\
&&\cdot(\sigma^{-1}(\alpha,X(\alpha-,t,x))\cdot(I+\nabla \sigma(\alpha,X(\alpha-,t,x))y)\cdot\nabla X(\alpha-,t,x)\cdot j_{\varepsilon}^{\flat}(y,r))^*\cdot h\\
&&-\zeta_{\varepsilon}^{1/2}(y)\cdot\sum_{i=1}^l(\sigma^{-1}(\alpha,X(\alpha-,t,x))\cdot\nabla \sigma(\alpha,X(\alpha-,t,x))e_i\cdot\nabla X(\alpha-,t,x)\cdot j_{\varepsilon}^{\flat}(y,r)\\
&&+\sigma^{-1}(\alpha,X(\alpha-,t,x))\cdot(I+\nabla \sigma(\alpha,X(\alpha-,t,x))y)\cdot\nabla X(\alpha-,t,x)\cdot\partial_i\zeta_{\varepsilon}^{1/2}(y)\cdot\xi(r))^*\cdot h]\rho(dr),
\de
where we have used:
\ce
\int_R\xi_i(r)\rho(dr)=0,\quad i=1,2,\cdots,l.
\de
By Proposition \ref{divergence}, we have
\ce
&&\delta_{\sharp}(F_{\varepsilon}^*(\tau)\cdot h)\\
&=&\int_t^{\tau}\int_{\cO}\varepsilon_{(\alpha,y)}^-(\widehat{\mE}[\delta_{\flat}(\varepsilon_{(\alpha,y,r)}^+(F_{\varepsilon}^*(\tau)\cdot h))(\alpha,y)])N(d\alpha,dy)\\
&=&\int_t^{\tau}\int_{\cO}\int_{R}\widehat{\mE}[-\theta_{\varepsilon}(y,r)\cdot(\sigma^{-1}(\alpha,X(\alpha-,t,x))\cdot(I+\nabla \sigma(\alpha,X(\alpha-,t,x))y)\\
&&\cdot\nabla X(\alpha-,t,x)\cdot j_{\varepsilon}^{\flat}(y,r))^*\cdot h\\
&&-\zeta_{\varepsilon}^{1/2}(y)\cdot\sum_{i=1}^l(\sigma^{-1}(\alpha,X(\alpha-,t,x))\cdot\nabla \sigma(\alpha,X(\alpha-,t,x))e_i\cdot\nabla X(\alpha-,t,x)\cdot j_{\varepsilon}^{\flat}(y,r)\\
&&+\sigma^{-1}(\alpha,X(\alpha-,t,x))\cdot(I+\nabla \sigma(\alpha,X(\alpha-,t,x))y)\cdot\nabla X(\alpha-,t,x)\cdot\partial_i\zeta_{\varepsilon}^{1/2}(y)\cdot\xi(r))^*\cdot h]\\
&&\rho(dr)N(d\alpha,dy),
\de
which completes the proof.
\end{proof}

\bl\label{gradient estimate 2}
Suppose that the assumptions {\bf (L)}, {\bf (FD)} and {\bf (FE)} are satisfied. Then for any $p\geq2$, there exists a positive constant $C$ such that for all $\tau\in[t,T]$ and $x,h\in\mR^d$,
\ce
\mE[|\widehat{\mE}[F_{\varepsilon}^*(\tau)\cdot h\cdot G_{\varepsilon}^{\sharp}(\tau)]|^p]\leq C(\tau-t)^p\varepsilon^{(5-\beta)p}|h|^p+C(\tau-t)\varepsilon^{5p-\beta}|h|^p
\de
and
\ce
\mE[|\delta_{\sharp}(F_{\varepsilon}^*(\tau)\cdot h)|^p]\leq C(\tau-t)\varepsilon^{2p-\beta}|h|^p+C(\tau-t)^p\varepsilon^{(2-\beta)p}|h|^p.
\de
\el

\begin{proof}

First, by Lemma \ref{bdg}, for all $t\in[0,T]$, we have 
\ce
&&\mE[|\nabla X(\tau,t,x)|^p]\\
&\leq&C+C\mE[\int_t^{\tau}\int_{\cO}|\nabla X(s-,t,x)|^2|y|^2\nu(dy)ds]^{p/2}+C\mE[\int_t^{\tau}\int_{\cO}|\nabla X(s-,t,x)|^p|y|^p\nu(dy)ds]\\
&&+C\mE[\int_t^{\tau}|\nabla X(\tau,t,x)|^pds]\\
&\leq&C+C\mE[\int_t^{\tau}|\nabla X(s,t,x)|^pds],
\de
from which we deduce by Gronwall's inequality that
\be\label{Jacobi matrix}
\sup_{\tau\in[t,T]}\sup_{x\in\mR^d}\mE[|\nabla X(\tau,t,x)|^p]\leq C.
\ee

We now prove the first conclusion. By Proposition \ref{marked Poisson measure}, we have
\ce
&&\widehat{\mE}[F_{\varepsilon}^*(\tau)\cdot h\cdot G_{\varepsilon}^{\sharp}(\tau)]\\
&=&\widehat{\mE}[(\int_t^{\tau}\int_{\cO}\int_R\sigma^{-1}(\alpha,X(\alpha-,t,x))\cdot(I+\nabla \sigma(\alpha,X(\alpha-,t,x))y)\\
&&\cdot\nabla X(\alpha-,t,x)\cdot j_{\varepsilon}^{\flat}(y,r)N\odot\rho(d\alpha,dy,dr))^*\cdot h\\
&&\cdot\int_t^{\tau}\int_{\cO}\int_R\zeta_{\varepsilon}^{1/2}(y)\cdot\sum_{i=1}^l\partial_i\zeta_{\varepsilon}(y)\cdot
\xi_i(r)N\odot\rho(d\alpha,dy,dr)]\\
&=&\int_t^{\tau}\int_{\cO}\int_R(\sigma^{-1}(\alpha,X(\alpha-,t,x))\cdot(I+\nabla \sigma(\alpha,X(\alpha-,t,x))y)\\
&&\cdot\nabla X(\alpha-,t,x)\cdot j_{\varepsilon}^{\flat}(y,r))^*\cdot h\cdot\zeta_{\varepsilon}^{1/2}(y)\cdot\sum_{i=1}^l\partial_i\zeta_{\varepsilon}(y)\cdot
\xi_i(r)\rho(dr)N(d\alpha,dy).
\de
Therefore using Kunita's inequality, (\ref{Jacobi matrix}), the assumption \textbf{(L)} and the assumption \textbf{(FE)}, we have
\ce
&&\mE[|\widehat{\mE}[F_{\varepsilon}^*(\tau)\cdot h\cdot G_{\varepsilon}^{\sharp}(\tau)]|^p]\\
&\leq&C\mE[|\int_t^{\tau}\int_{\cO}|\int_R(\sigma^{-1}(\alpha,X(\alpha-,t,x))\cdot(I+\nabla \sigma(\alpha,X(\alpha-,t,x))y)\\
&&\cdot\nabla X(\alpha-,t,x)\cdot j_{\varepsilon}^{\flat}(y,r))^*\cdot h\cdot\zeta_{\varepsilon}^{1/2}(y)\cdot\sum_{i=1}^l\partial_i\zeta_{\varepsilon}(y)\cdot
\xi(r)\rho(dr)|\nu(dy)d\alpha|^p]\\
&&+C\mE[\int_t^{\tau}\int_{\cO}|\int_R(\sigma^{-1}(\alpha,X(\alpha-,t,x))\cdot(I+\nabla \sigma(\alpha,X(\alpha-,t,x))y)\\
&&\cdot\nabla X(\alpha-,t,x)\cdot j_{\varepsilon}^{\flat}(y,r))^*\cdot h\cdot\zeta_{\varepsilon}^{1/2}(y)\cdot\sum_{i=1}^l\partial_i\zeta_{\varepsilon}(y)\cdot
\xi(r)\rho(dr)|^p\nu(dy)d\alpha]\\
&\leq&C(\int_t^{\tau}\int_{\{|u|\leq\varepsilon\}}|u|^5\nu(du)ds)^p|h|^p\\
&&+C\int_t^{\tau}\int_{\{|u|\leq\varepsilon\}}|u|^{5p}\nu(du)ds|h|^p\\
&\leq&C(\tau-t)^p\varepsilon^{(5-\beta)p}|h|^p+C(\tau-t)\varepsilon^{5p-\beta}|h|^p.
\de

Now we start to prove the second conclusion. Using the fact that 
\ce
\int_R\xi_i(r)\rho(dr)=0,\quad i=1,2,\cdots,l,
\de
and applying Lemma \ref{divergence operator}, Kunita's inequality, the assumption \textbf{(L)} and the assumption \textbf{(FD)}, we obtain
\ce
&&\mE[|\delta_{\sharp}(F_{\varepsilon}^*(\tau)\cdot h)|^p]\\
&\leq&C\mE[\int_t^{\tau}\int_{\cO}|\int_{R}\widehat{\mE}[\theta_{\varepsilon}(y,r)\cdot(\sigma^{-1}(\alpha,X(\alpha-,t,x))\cdot(I+\nabla \sigma(\alpha,X(\alpha-,t,x))y)\\
&&\cdot\nabla X(\alpha-,t,x)\cdot j_{\varepsilon}^{\flat}(y,r))^*\cdot h]\rho(dr)|^p\nu(dy)d\alpha]\\
&&+C\mE[|\int_t^{\tau}\int_{\cO}\int_{R}\widehat{\mE}[\theta_{\varepsilon}(y,r)\cdot(\sigma^{-1}(\alpha,X(\alpha-,t,x))\cdot(I+\nabla \sigma(\alpha,X(\alpha-,t,x))y)\\
&&\cdot\nabla X(\alpha-,t,x)\cdot j_{\varepsilon}^{\flat}(y,r))^*\cdot h]\rho(dr)\nu(dy)d\alpha|^p]\\
&\leq&C\mE[\int_t^{\tau}\int_{\cO}|\int_{R}\widehat{\mE}[\zeta_{\varepsilon}^{1/2}(y)\cdot\sum_{i=1}^l\\
&&(\sigma^{-1}(\alpha,X(\alpha-,t,x))\cdot\nabla \sigma(\alpha,X(\alpha-,x))e_i\cdot\nabla X(\alpha-,t,x)\cdot j_{\varepsilon}^{\flat}(y,r)\\
&&+\sigma^{-1}(\alpha,X(\alpha-,t,x))\cdot(I+\nabla \sigma(\alpha,X(\alpha-,t,x))y)\cdot\nabla X(\alpha-,t,x)\cdot\partial_i\zeta_{\varepsilon}^{1/2}(y)\cdot\xi(r))^*\\
&&\cdot h]\rho(dr)|^p\nu(dy)d\alpha]\\
&&+C\mE[|\int_t^{\tau}\int_{\cO}\int_{R^2}\widehat{\mE}[\zeta_{\varepsilon}^{1/2}(y)\cdot\sum_{i=1}^l\\
&&(\sigma^{-1}(\alpha,X(\alpha-,t,x))\cdot\nabla \sigma(\alpha,X(\alpha-,t,x))e_i\cdot\nabla X(\alpha-,t,x)\cdot j_{\varepsilon}^{\flat}(y,r)\\
&&+\sigma^{-1}(\alpha,X(\alpha-,t,x))\cdot(I+\nabla \sigma(\alpha,X(\alpha-,t,x))y)\cdot\nabla X(\alpha-,t,x)\cdot\partial_i\zeta_{\varepsilon}^{1/2}(y)\cdot\xi(r))^*\\
&&\cdot h]\rho(dr)\nu(dy)d\alpha|^p]\\
&\leq&C\int_t^{\tau}\int_{\{|u|\leq\varepsilon\}}(|u|^3+|u^2|)^p\nu(du)ds|h|^p\\
&&+C(\int_t^{\tau}\int_{\{|u|\leq\varepsilon\}}(|u|^3+|u^2|)\nu(du)ds)^p|h|^p\\
&\leq&C(\tau-t)\varepsilon^{2p-\beta}|h|^p+C(\tau-t)^p\varepsilon^{(2-\beta)p}|h|^p,
\de
and this completes the proof.
\end{proof}

Let
\ce
U_{\tau,\varepsilon}^h=\delta_{\sharp}(F_{\varepsilon}^*(\tau)\cdot h)\cdot G_{\varepsilon}^{-1}(\tau)-\widehat{\mE}[F_{\varepsilon}^*(\tau)\cdot h\cdot (G_{\varepsilon}^{-1}(\tau))^{\sharp}]
\de
and
\ce
U_{\tau}^h:=U_{\tau,(\tau-t)^{\frac{1}{\beta}}}^h.
\de
Therefore using the assumption (\textbf{L}-c) and Lemma \ref{gradient estimate 2}, one can prove the  following result.

\bl\label{Bismut estimate}
Suppose that the assumptions {\bf (L)}, {\bf (FD)} and {\bf (FE)} are satisfied. Then for any $\tau\in[t,T]$ and $h\in\mR^d$, we have
\ce
&&\mE[|U_{\tau,\epsilon}^h|^p]^{1/p}\\
&\leq& C((\tau-t)\varepsilon^{4p-\beta}+(\tau-t)^{2p}\varepsilon^{(2-\beta)2p})^{1/2p}(((\tau-t)\varepsilon^{3-\beta})^{-2p}+(\tau-t)^{-\frac{6p}{\beta}})^{1/2p}|h|\\
&&+C((\tau-t)^{2p}\varepsilon^{(5-\beta)2p}+(\tau-t)\varepsilon^{10p-\beta})^{1/2p}(((\tau-t)\varepsilon^{3-\beta})^{-4p}+(\tau-t)^{-\frac{12p}{\beta}})^{1/2p}|h|,
\de
and in particular,
\ce
\mE[|U_{\tau}^h|^p]^{1/p}\leq C(\tau-t)^{-1/\beta}|h|.
\de
\el

\subsection{Proof of Theorem \ref{Bismut formula}}

For the sake of simplicity, the following notations will be used:
\ce
(X(\tau,x),Y(\tau,x),Z(\tau,x,u))=(X(\tau,t,x),Y(\tau,t,x),X(\tau,t,x,u)).
\de
Let $\nabla_xY$ and $\nabla_xZ$ be the partial derivatives with respect to $x$. First of all, using Lemma \ref{gradient of stochastic integration}, for any $0\leq t\leq s\leq T$, we have
\ce
&&(\psi(s,X(s,x),Y(s,x),\int_{\cO}Z(s,x,u)l(u)\nu(du)))^{\sharp}\\
&=&\nabla_x\psi(s,X(s,x),Y(s,x),\int_{\cO}Z(s,x,u)l(u)\nu(du))(X(s,x))^{\sharp}\\
&&+\nabla_y\psi(s,X(s,x),Y(s,x),\int_{\cO}Z(s,x,u)l(u)\nu(du))(Y(s,x))^{\sharp}\\
&&+\nabla_z\psi(s,X(s,x),Y(s,x),\int_{\cO}Z(s,x,u)l(u)\nu(du))\int_{\cO}(Z(s,x,u))^{\sharp}l(u)\nu(du).
\de
This implies by Lemma \ref{derivative initial} that
\be\label{initial Malliavin}
\begin{split}
&\widehat{\mE}[(\psi(s,X(s,x),Y(s,x),\int_{\cO}Z(s,x,u)l(u)\nu(du)))^{\sharp}F_{\varepsilon}^*(s)G_{\varepsilon}^{-1}(s)]\\
=&\widehat{\mE}[\nabla_x\psi(s,X(s,x),Y(s,x),\int_{\cO}Z(s,x,u)l(u)\nu(du))(X(s,x))^{\sharp}F_{\varepsilon}^*(s)G_{\varepsilon}^{-1}(s)\\
&+\nabla_y\psi(s,X(s,x),Y(s,x),\int_{\cO}Z(s,x,u)l(u)\nu(du))(Y(s,x))^{\sharp}F_{\varepsilon}^*(s)G_{\varepsilon}^{-1}(s)\\
&+\nabla_z\psi(s,X(s,x),Y(s,x),\int_{\cO}Z(s,x,u)l(u)\nu(du))\\
&\int_{\cO}(Z(s,x,u))^{\sharp}l(u)\nu(du)F_{\varepsilon}^*(s)G_{\varepsilon}^{-1}(s)]\\
=&\nabla_x\psi(s,X(s,x),Y(s,x),\int_{\cO}Z(s,x,u)l(u)\nu(du))\nabla_xX(s,x)\\
&+\nabla_y\psi(s,X(s,x),Y(s,x),\int_{\cO}Z(s,x,u)l(u)\nu(du))\nabla_xY(s,x)\\
&+\nabla_z\psi(s,X(s,x),Y(s,x),\int_{\cO}Z(s,x,u)l(u)\nu(du))\int_{\cO}\nabla_xZ(s,x,u)l(u)\nu(du).
\end{split}
\ee
On the other hand, it is easy to see that the directional derivative process in the direction $h\in\mR^d$, $\{\nabla_xY(\tau,x)h,\nabla_xZ(\tau,x,u)h;\tau\in[0,T],x,h\in\mR^d,u\in\cO\}$, solves the following backward SDEs with jumps:
\be\label{initial}
\begin{split}
&\nabla_xY(\tau,x)h+\int_{\tau}^T\int_{\cO}\nabla_xZ(s-,x,u)h\widetilde{N}(ds,du)\\
=&-\int_{\tau}^T\nabla_x\psi(s,X(s,x),Y(s,x),\int_{\cO}Z(s,x,u)l(u)\nu(du))\nabla_xX(s,x)hds\\
&-\int_{\tau}^T\nabla_y\psi(s,X(s,x),Y(s,x),\int_{\cO}Z(s,x,u)l(u)\nu(du))\nabla_xY(s,x)hds\\
&-\int_{\tau}^T\nabla_z\psi(s,X(s,x),Y(s,x),\int_{\cO}Z(s,x,u)l(u)\nu(du))\\
&\int_{\cO}\nabla_xZ(s,x,u)l(u)\nu(du)hds\\
&+\nabla\phi(X(T,x))\nabla_xX(T,x)h.
\end{split}
\ee
Hence taking expectations on both sides in (\ref{initial}), and using (\ref{initial Malliavin}) and Lemma \ref{derivative initial}, we have
\ce
&&\mE[\nabla_xY(\tau,x)h]\\
&=&-\mE[\int_{\tau}^T\nabla_x\psi(s,X(s,x),Y(s,x),\int_{\cO}Z(s,x,u)l(u)\nu(du))\nabla_xX(s,x)hds\\
&&+\int_{\tau}^T\nabla_y\psi(s,X(s,x),Y(s,x),\int_{\cO}Z(s,x,u)l(u)\nu(du))\nabla_xY(s,x)hds\\
&&+\int_{\tau}^T\nabla_z\psi(s,X(s,x),Y(s,x),\int_{\cO}Z(s,x,u)l(u)\nu(du))\int_{\cO}\nabla_xZ(s,x,u)l(u)\nu(du)hds]\\
&&+\mE[\nabla\phi(X(T,x))\nabla_xX(T,x)h]\\
&=&-\mE\widehat{\mE}[\int_{\tau}^T(\psi(s,X(s,x),Y(s,x),\int_{\cO}Z(s,x,u)l(u)\nu(du)))^{\sharp}\cdot F_{\varepsilon}^*(s)\cdot h\cdot G_{\varepsilon}^{-1}(s)ds]\\
&&+\mE\widehat{\mE}[\nabla\phi(X(T,x))(X(T,x))^{\sharp}\cdot F_{\varepsilon}^*(T)\cdot h\cdot G_{\varepsilon}^{-1}(T)]\\
&=&-\mE\widehat{\mE}[\int_{\tau}^T(\psi(s,X(s,x),Y(s,x),\int_{\cO}Z(s,x,u)l(u)\nu(du)))^{\sharp}\cdot F_{\varepsilon}^*(s)\cdot h\cdot G_{\varepsilon}^{-1}(s)ds]\\
&&+\mE\widehat{\mE}[(\phi(X(T,x)))^{\sharp}\cdot F_{\varepsilon}^*(T)\cdot h\cdot G_{\varepsilon}^{-1}(T)]\\
&=&-\mE[\int_{\tau}^T\psi(s,X(s,x),Y(s,x),\int_{\cO}Z(s,x,u)l(u)\nu(du))\cdot[\delta_{\sharp}(F_{\varepsilon}^*(s)\cdot h)\cdot G_{\varepsilon}^{-1}(s)\\
&&-\widehat{\mE}[F_{\varepsilon}^*(s)\cdot h\cdot (G_{\varepsilon}^{-1}(s))^{\sharp}]]ds]\\
&&+\mE[\phi(X(t,x))\cdot[\delta_{\sharp}(F_{\varepsilon}^*(T)\cdot h)\cdot G_{\varepsilon}^{-1}(T)-\widehat{\mE}[F_{\varepsilon}^*(T)\cdot h\cdot (G_{\varepsilon}^{-1}(T))^{\sharp}]]].
\de
Then the proof is complete.

\subsection{Proof of Theorem \ref{gradient estimate}}

We are now ready to prove Theorem \ref{gradient estimate}. For simplicity the desired gradient estimate will be proved under the additional assumption that $t$ takes values in an interval $[T-\delta,T]$ sufficiently small. By a standard technique (introducing an exponential weight) the same argument gives the result in the whole $[0,T ]$.

Since  $Y(t,t,x)$ is deterministic, so is $\nabla_xY(t,t,x)$. Using Theorem \ref{BSDE existence}, the assumption {\bf (BL)}, and Theorem \ref{FSDE existence}, we have
\be\label{BSDE estimate}
\begin{split}
|Y(t,t,x)|^2=&\mE[|Y(t,t,x)|^2]\\
\leq&\mE[\sup_{\tau\in[0,T]}|Y(\tau,t,x)|^2]\\
\leq&C\mE[\int_0^T|\psi(s,X(s,t,x),0,0)|^2ds]+C\mE[|\phi(X(T,t,x))|^2]\\
\leq&C(1+|x|)^{\mu}.
\end{split}
\ee

Now using Theorem \ref{Bismut formula} and Theorem \ref{PDE existence BD}, it holds that
\ce
&&\mE[\nabla Y(\tau,t,x)h]\\
&=&-\mE[\int_{\tau}^T\psi(s,X(s,t,x),Y(s,t,x),\int_{\cO}Z(s,t,x,u)l(u)\nu(du))\\
&&\cdot[\delta_{\sharp}(F_{(T-t)^{\frac{1}{\beta}}}^*(s)\cdot h)\cdot G_{(T-t)^{\frac{1}{\beta}}}^{-1}(s)-\widehat{\mE}[F_{(T-t)^{\frac{1}{\beta}}}^*(s)\cdot h\cdot (G_{(T-t)^{\frac{1}{\beta}}}^{-1}(s))^{\sharp}]]ds]\\
&&+\mE[\phi(X(T,t,x))\cdot[\delta_{\sharp}(F_{(T-t)^{\frac{1}{\beta}}}^*(T)\cdot h)\cdot G_{(T-t)^{\frac{1}{\beta}}}^{-1}(T)\\
&&-\widehat{\mE}[F_{(T-t)^{\frac{1}{\beta}}}^*(T)\cdot h\cdot (G_{(T-t)^{\frac{1}{\beta}}}^{-1}(T))^{\sharp}]]]\\
&=&-\mE[\int_{\tau}^T\psi(s,X(s,t,x),Y(s,t,x),\int_{\cO}(v(s,X(s,t,x)+\sigma(s,X(s,t,x))u)\\
&&-v(s,X(s,t,x)))l(u)\nu(du))\cdot[\delta_{\sharp}(F_{(T-t)^{\frac{1}{\beta}}}^*(s)\cdot h)\cdot G_{(T-t)^{\frac{1}{\beta}}}^{-1}(s)\\
&&-\widehat{\mE}[F_{(T-t)^{\frac{1}{\beta}}}^*(s)\cdot h\cdot (G_{(T-t)^{\frac{1}{\beta}}}^{-1}(s))^{\sharp}]]ds]\\
&&+\mE[\phi(X(T,t,x))\cdot[\delta_{\sharp}(F_{(T-t)^{\frac{1}{\beta}}}^*(T)\cdot h)\cdot G_{(T-t)^{\frac{1}{\beta}}}^{-1}(T)\\
&&-\widehat{\mE}[F_{(T-t)^{\frac{1}{\beta}}}^*(T)\cdot h\cdot (G_{(T-t)^{\frac{1}{\beta}}}^{-1}(T))^{\sharp}]]]\\
&=&-\mE[\int_{\tau}^T\psi(s,X(s,t,x),Y(s,t,x),\int_{\cO}(v(s,X(s,t,x)+\sigma(s,X(s,t,x))u)\\
&&-v(s,X(s,t,x)))l(u)\nu(du))\cdot U_s^hds]\\
&&+\mE[\phi(X(T,t,x))\cdot U_T^h],
\de
where $v(t,x)=Y(t,t,x)$. Then it follows from the assumption {\bf (BL)} that
\ce
&&|\nabla Y(t,t,x)h|\\
&\leq&\mE[|\int_t^T\psi(s,X(s,t,x),0,0)\cdot U_s^hds|]\\
&&+\mE[|\int_t^T(\psi(s,X(s,t,x),Y(s,t,x),\int_{\cO}(v(s,X(s,t,x)+\sigma(s,X(s,t,x))u)\\
&&-v(s,X(s,t,x)))l(u)\nu(du))\\
&&-\psi(s,X(s,t,x),0,0))\cdot U_s^hds|]\\
&&+\mE[\phi(X(T,t,x))\cdot U_T^h]\\
&\leq&C_{BL}\mE[\int_t^T(1+|X(s,t,x)|)^{\mu}|U_s^h|ds]\\
&&+C_{BL}\mE[\int_t^T|Y(s,t,x)||U_s^h|ds]\\
&&+C_{BL}\mE[\int_t^T|\int_{\cO}(v(s,X(s,t,x)+\sigma(s,X(s,t,x))u)\\
&&-v(s,X(s,t,x)))l(u)\nu(du)||U_s^h|ds]\\
&&+C_{BL}\mE[(1+|X(T,t,x)|)^{\mu}|U_T^h|]\\
&=:&I_1(t,x)+I_2(t,x)+I_3(t,x)+I_4(t,x).
\de
Define
\ce
|||\nabla_xY|||=\sup_{t\in[T-\delta,T]}\sup_{x\in\mR^d}(T-t)^{1/\beta}(1+|x|)^{-\mu}|\nabla_xY(t,t,x)|,
\de
and $|||I_k|||$, $k=1,2,3,4$, is defined similarly.

For the term $I_4$, using H\"{o}lder's inequality, Theorem \ref{FSDE existence} and Lemma \ref{Bismut estimate}, we have
\ce
I_4(t,x)&\leq&C_{BL}\mE[(1+X(T,t,x))^{2\mu}]^{1/2}\mE[|U_T^h|^2]^{1/2}\\
&\leq&C(1+|x|)^{\mu}(T-t)^{-1/\beta}|h|,
\de
and hence $|||I_4|||<\infty$.

Now we deal with the term $I_3$. In view of H\"{o}lder's inequality, the mean value theorem, the assumption {\bf (FE)}, Lemma \ref{Bismut estimate} and Theorme \ref{FSDE existence}, it holds that
\ce
&&I_3(t,x)\\
&\leq&C_{BL}\int_t^T\mE[|\int_{\cO}(v(s,X(s,t,x)+\sigma(s,X(s,t,x))u)\\
&&-v(s,X(s,t,x)))l(u)\nu(du)|^2]^{1/2}\mE[|U_s^h|^2]^{1/2}ds\\
&\leq&C_{BL}\int_t^T\mE[|\int_{\cO}|\nabla v(s,X(s,t,x)+(1-\theta)\sigma(s,X(s,t,x))u)\sigma(s,X(s,t,x))|\\
&&|ul(u)|\nu(du)|^2]^{1/2}\mE[|U_s^h|^2]^{1/2}ds\\
&\leq&C_{BL}C_{FB}\int_t^T\mE[\int_{\cO}|\nabla v(s,X(s,t,x)+(1-\theta)\sigma(s,X(s,t,x))u)|\\
&&|ul(u)|\nu(du)|^2]^{1/2}\mE[|U_s^h|^2]^{1/2}ds\\
&\leq&C\int_{\cO}|ul(u)|\nu(du)\\
&&\int_t^T(T-s)^{-1/\beta}\mE[(1+|X(s,t,x)|+(1-\theta)|\sigma(s,X(s,t,x))|)^{2\mu}]^{1/2}\\
&&(s-t)^{-1/\beta}ds|h||||\nabla_xv|||\\
&\leq&C(1+|x|)^{\mu}\int_t^T(T-s)^{-1/\beta}(s-t)^{-1/\beta}ds|h|||\nabla_xY|||\\
&=&C(1+|x|)^{\mu}(T-t)^{-\frac{2}{\beta}+1}\int_0^1(1-z)^{-\frac{1}{\beta}}z^{-\frac{1}{\beta}}dz|h|||\nabla_xY|||,
\de
where $0<\theta<1$. This implies that
\ce
|||I_3|||\leq C\delta^{-\frac{1}{\beta}+1}|||\nabla_xY|||.
\de

Similarly, using Theorem \ref{BSDE existence}, we have
\ce
|||I_1|||\leq C\delta,\quad|||I_2|||\leq C\delta.
\de
Combining all the above  estimates, we obtain
\ce
|||\nabla_xY|||\leq C+C\delta+C\delta^{-\frac{1}{\beta}+1}|||\nabla_xY|||.
\de
Since $\beta\in(1,2)$, we can take $\delta$ small enough such that
\ce
|||\nabla_xY|||\leq C.
\de
This is the desired result, and hence we concludes the proof.

\section{Nonlocal Quasi-Linear Integral-PDEs}

\subsection{Notations and Notions}

We introduce the semigroup $\{P_{t,\tau};0\leq t\leq\tau\leq T\}$, acting on the space of Borel functions $\phi:\mR^d\rightarrow\mR$ having polynomial growth, defined by:
\ce
P_{t,\tau}[\phi](x):=\mE[\phi(X(\tau,t,x))],\quad x\in\mR^d.
\de
By $\cL_t$ we denote the formal generator of $P_{t,\tau}$, namely:
\ce
\cL_t[\phi](x):=b(t,x)\cdot\nabla\phi(x)+p.v.\int_{\cO}(\phi(x+\sigma(t,x)u)-\phi(x))\nu(du),
\de
where p.v. stands for the Cauchy principle value.

We consider the following nonlocal quasi-linear integral-PDEs: for any $t\in[0,T]$ and $x\in\mR^d$,
\be\label{PDE}
\begin{cases}
\frac{\partial v(t,x)}{\partial t}+\cL_t[v(t,\cdot)](x)\\
\quad\quad\quad=\psi(t,x,v(t,x),\int_\cO(v(t,x+\sigma(t,x)u)-v(t,x))l(u)\nu(du)),\\
v(T,x)=\phi(x).
\end{cases}
\ee

Now we give the definition of the mild solution to Eq. (\ref{PDE}).

\bd
We say that a continuous function $v:[0,T]\times\mR^d\rightarrow\mR$ is a mild solution of Eq. (\ref{PDE}) if the following conditions hold:
\begin{enumerate}
\item
$v\in\cC^{0,1}([0,T]\times\mR^d;\mR)$;
\item
for all $t\in[0,T]$ and $x\in\mR^d$, we have
\ce
|v(t,x)|\leq C(1+|x|)^C,\quad|\nabla_xv(t,x)|\leq Cf(t)(1+|x|)^C,
\de
for some constant $C>0$, and some real function $f$ satisfying $\int_0^Tf(t)dt<\infty$;
\item
for any $t\in[0,T]$ and $x\in\mR^d$, the following equality holds:
\ce
v(t,x)=-\int_t^TP_{t,\tau}[\psi(\tau,\cdot,v(\tau,\cdot),\int_{\cO}(v(\tau,\cdot+\sigma(\tau,\cdot)u)-v(\tau,\cdot))l(u)du)](x)d\tau+P_{t,T}[\phi](x).
\de
\end{enumerate}
\ed

\subsection{Joint Quadratic Variation}

Since we are considering mild solutions here, in order to prove the uniqueness of the solution to Eq. (\ref{PDE}), a result on  joint quadratic variation should be introduced. The definition of generalized joint quadratic variation we consider in the present paper has been first introduced in \cite{RV}, and it is shown in \cite[Proposition 1.1]{RV} that the standard definition of joint quadratic variation coincides with the generalized joint quadratic variation defined below.

\bd
Given a couple of $\mR$-valued càdlàg processes $(X(t),Y(t))$, $t\geq0$, we define their joint quadratic variation on $[0,T]$, to be
\ce
\<X(t),Y(t)\>_{[0,T]}:=\mP-\lim_{\epsilon\downarrow0}\frac{1}{\epsilon}\int_0^T(X((t+\epsilon)\wedge T)-X(t))(Y((t+\epsilon)\wedge T)-Y(t))dt,
\de
where $\mP-\lim$ denotes the limit to be taken in probability.
\ed

The following lemma can be found in \cite[Theorem 3.2]{CPO}.

\bl\label{quadratic variation}
Assume that $v:[0,T]\times\mR^d\rightarrow\mR$ is locally Lipschitz with respect to the second variable and with at most polynomial growth, i.e., there exist $C>0$ and $m\geq0$, such that, for any $t\in[0,T]$ and $x,y\in\mR^d$, $v$ satisfies
\ce
|v(t,x)-v(t,y)|&\leq& C|x-y|(1+|x|+|y|)^m,\\
|v(t,0,0)|&\leq& C.
\de
Then, for every $t\in[0,T]$ and $x\in\mR^d$, the process
\ce
\{v(s,X(s,t,x));s\in[t,T]\}
\de
admits a joint quadratic variation on the interval $[t,T]$ with
\ce
J(s)=\int_t^s\int_{\cO}u\widetilde{N}(ds,du)
\de
given by
\ce
&&\<v(\cdot,X(\cdot,t,x)),J(\cdot)\>_{[t,T]}\\
&=&\int_t^T\int_{\cO}(v(\tau,X(\tau,t,x)+\sigma(\tau,X(\tau,t,x))u)-v(\tau,X(\tau,t,x)))uN(d\tau,du).
\de
\el

\subsection{Nonlocal Quasi-Linear Integral-PDEs with the Assumption {\bf (BD)}}

\bt\label{PDE existence BD}
Suppose that $\beta\in(1,2)$ and the assumptions {\bf (L)}, {\bf (FD)}, {\bf (FE)}, {\bf (BL)} and {\bf (BD)} are satisfied. Then there exists a unique mild solution $v$ of Eq. (\ref{PDE}) given by the formula
\ce
v(t,x)=Y(t,t,x),
\de
and the following representation formulas are satisfied:
\ce
Y(\tau,t,x)&=&v(\tau,X(\tau,t,x)),\quad\mP-\text{a.s.},\\
Z(\tau,t,x,u)&=&v(\tau,X(\tau,t,x)+\sigma(\tau,X(\tau,t,x))u)-v(\tau,X(\tau,t,x)),\quad d\tau d\nu d\mP-\text{a.e.}
\de
\et

This proof is similar to the proof of Theorem 4.8 in \cite{CPO}. For the convenience of readers, we include the proof here.

\begin{proof}[Proof of Existence]

First of all, the representation of $Y$ follows from a standard consequence of uniqueness of the solution of Eq. (\ref{BSDE}). For the representation of $Z$, we argue as follows. Using the standard definition of joint quadratic variation, we have
\be\label{variation Y}
\<Y(\cdot,t,x),J(\cdot)\>_{[t,T]}=\int_t^T\int_{\cO}Z(\tau,t,x,u)uN(d\tau,du).
\ee
On the other hand, it can be seen from Lemma \ref{quadratic variation} that
\be\label{variation u}
\begin{split}
&\<v(\cdot,X(\cdot,t,x)),J(\cdot)\>_{[t,T]}\\
=&\int_t^T\int_{\cO}v(\tau,X(\tau,t,x)+\sigma(\tau,X(\tau,t,x))u)-v(\tau,X(\tau,t,x))uN(d\tau,du).
\end{split}
\ee
Comparing Eq. (\ref{variation Y}) and Eq. (\ref{variation u}), the representation of $Z$ is obtained.

Now using the fact that $v(t,x)=Y(t,t,x)$ is deterministic, we have
\ce
&&v(t,x)\\
&=&-\mE[\int_t^T\psi(\tau,X(\tau,t,x),v(\tau,X(\tau,t,x)),\int_{\cO}Z(\tau,x,t,u)l(u)\nu(du))d\tau]\\
&&+\mE[\phi(X(T,t,x))]\\
&=&-\mE[\int_t^T\psi(\tau,X(\tau,t,x),v(\tau,X(\tau,t,x)),\\
&&\int_{\cO}(v(\tau,X(\tau,t,x)+\sigma(\tau,X(\tau,t,x))u)-v(\tau,X(\tau,t,x)))l(u)\nu(du))d\tau]\\
&&+\mE[\phi(X(T,t,x))]\\
&=&-\int_t^TP_{t,\tau}[\psi(\tau,\cdot,v(\tau,\cdot),\int_{\cO}(v(\tau,\cdot+\sigma(\tau,\cdot)u)-v(\tau,\cdot))l(u)\nu(du))](x)d\tau\\
&&+P_{t,T}[\phi](x).
\de
It can be seen from (\ref{BSDE estimate}) that
\ce
|v(t,x)|\leq C(1+|x|^{\mu}),\quad\forall t\in[0,T],~\forall x\in\mR^d.
\de
Applying the estimate in Theorem \ref{gradient estimate}, we also obtain
\ce
|\nabla_xv(t,x)|\leq C(T-t)^{-1/\beta}(1+|x|)^{\mu},\quad\forall t\in[0,T],~\forall x\in\mR^d.
\de
Then the proof of existence of Eq. (\ref{PDE}) is complete.

\end{proof}

\begin{proof}[Proof of Uniqueness]

Let $v$ be a mild solution of Eq. (\ref{PDE}). Using the definition of mild solution and the definition of $P_{t,\tau}$, for every $s\in[t,T]$ and $x\in\mR^d$,
\ce
v(s,x)&=&-\mE[\int_s^T\psi(\tau,X(\tau,s,x),v(\tau,X(\tau,s,x)),\\
&&\int_{\cO}(v(\tau,X(\tau,s,x)+\sigma(\tau,X(\tau,s,x))u)-v(\tau,X(\tau,s,x)))l(u)\nu(du))d\tau]\\
&&+\mE[\phi(X(T,s,x))].
\de
Since $X(\tau,s,x)$ is $\mathcal{F}_s$-measurable, we have
\ce
v(s,x)&=&-\mE[\int_s^T\psi(\tau,X(\tau,s,x),v(\tau,X(\tau,s,x)),\\
&&\int_{\cO}(v(\tau,X(\tau,s,x)+\sigma(X(\tau,s,x))u)-v(\tau,X(\tau,s,x)))l(u)\nu(du))d\tau|\mathcal{F}_s]\\
&&+\mE[\phi(X(T,s,x))|\mathcal{F}_s].
\de
Now replacing $x$ by $X(s,t,x)$ and using the fact that
\ce
X(\tau,s,X(s,t,x))=X(\tau,t,x),\quad\mP-\text{a.s. for}~\tau\in[s,T],
\de
we obtain
\ce
v(s,X(s,t,x))&=&-\mE[\int_s^T\psi(\tau,X(\tau,t,x),v(\tau,X(\tau,t,x)),\\
&&\int_{\cO}(v(\tau,X(\tau,t,x)+\sigma(\tau,X(\tau,t,x))u)-v(\tau,X(\tau,t,x)))l(u)\nu(du))d\tau|\mathcal{F}_s]\\
&&+\mE[\phi(X(T,t,x))|\mathcal{F}_s]\\
&=&\int_t^s\psi(\tau,X(\tau,t,x),v(\tau,X(\tau,t,x)),\\
&&\int_{\cO}(v(\tau,X(\tau,t,x)+\sigma(\tau,X(\tau,t,x))u)-v(\tau,X(\tau,t,x)))l(u)\nu(du))d\tau\\
&&+\mE[\xi|\mathcal{F}_s],
\de
where
\ce
\xi&=&-\int_t^T\psi(\tau,X(\tau,t,x),v(\tau,X(\tau,t,x)),\\
&&\int_{\cO}(v(\tau,X(\tau,t,x)+\sigma(\tau,X(\tau,t,x))u)-v(\tau,X(\tau,t,x)))l(u)\nu(du))d\tau\\
&&+\phi(X(T,t,x)).
\de
Noticing that $\mE[\phi(X(T,t,x))|\mathcal{F}_s]=v(t,x)$ and using the fact that $\xi\in L^2(\Omega;\mR)$ is $\mathcal{F}_T$-measurable and applying the martingale representation theorem, there exists a predictable process $U:\Omega\times[0,T]\times\mR\rightarrow\mR$ satisfying
\ce
\mE[\int_0^T\int_{\cO}|U(t,u)|^2\nu(du)dt]<\infty,
\de
such that
\ce
\mE[\xi|\mathcal{F}_s]=\int_t^s\int_{\cO}U(\tau,u)\widetilde{N}(d\tau,du)+v(t,x).
\de
This implies that the process $\{v(s,X(s,t,x));s\in[t,T],x\in\mR^d\}$ is a real continuous semimartingale with canonical decomposition:
\be\label{decomposition}
\begin{split}
&v(s,X(s,t,x))\\
=&\int_t^s\psi(\tau,X(\tau,t,x),v(\tau,X(\tau,t,x)),\\
&\int_{\cO}(v(\tau,X(\tau,t,x)+\sigma(\tau,X(\tau,t,x))u)-v(\tau,X(\tau,t,x)))l(u)\nu(du))d\tau\\
&+\int_t^s\int_{\cO}U(\tau,u)\widetilde{N}(d\tau,du)+v(t,x).
\end{split}
\ee
Now using Lemma \ref{quadratic variation}, it holds that
\ce
&&\int_t^s\int_{\cO}(v(\tau,X(\tau,t,x)+\sigma(\tau,X(\tau,t,x))u)-v(\tau,X(\tau,t,x)))uN(d\tau,du)\\
&=&\int_t^s\int_{\cO}U(\tau,u)uN(d\tau,du).
\de
Therefore, for $d\tau d\nu$-a.e. $(\tau,u)\in[t,T]\times\cO$, we have, $\mP$-a.s.,
\ce
(v(\tau,X(\tau,t,x)+\sigma(\tau,X(\tau,t,x))u)-v(\tau,X(\tau,t,x)))=U(\tau,u).
\de
Substituting this into (\ref{decomposition}), we obtain that for any $s\in[t,T]$,
\ce
&&v(s,X(s,t,x))\\
&=&\int_t^s\psi(\tau,X(\tau,t,x),v(\tau,X(\tau,t,x)),\\
&&\int_{\cO}(v(\tau,X(\tau,t,x)+\sigma(\tau,X(\tau,t,x))u)-v(\tau,X(\tau,t,x)))l(u)\nu(du))d\tau\\
&&+\int_t^s\int_{\cO}(v(\tau,X(\tau,t,x)+\sigma(\tau,X(\tau,t,x))u)-v(\tau,X(\tau,t,x)))\widetilde{N}(d\tau,du)+v(t,x).
\de
This implies by $v(T,X(T,t,x))=\phi(X(T,t,x))$ that for any $s\in[t,T]$,
\ce
&&v(s,X(s,t,x))+\int_s^T\int_{\cO}(v(\tau,X(\tau,t,x)+\sigma(\tau,X(\tau,t,x))u)-v(\tau,X(\tau,t,x)))\widetilde{N}(d\tau,du)\\
&=&\phi(X(T,t,x))-\int_s^T\psi(\tau,X(\tau,t,x),v(\tau,X(\tau,t,x)),\\
&&\int_{\cO}(v(\tau,X(\tau,t,x)+\sigma(\tau,X(\tau,t,x))u)-v(\tau,X(\tau,t,x)))l(u)\nu(du))d\tau.
\de
Comparing with Eq. (\ref{BSDE}), it is easy to see that the pair
\ce
(Y(s,t,x),Z(s,t,x,u))
\de
and
\ce
(v(s,X(s,t,x)),v(s,X(s,t,x)+\sigma(s,X(s,t,x))u)-v(s,X(s,t,x)))
\de
solve the same equation. By uniqueness, we have
\ce
Y(s,t,x)=v(s,X(s,t,x))
\de
and
\ce
Z(s,t,x,u)=v(s,X(s,t,x)+\sigma(s,X(s,t,x))u)-v(s,X(s,t,x)).
\de
Then we complete the proof.

\end{proof}

\subsection{Nonlocal Quasi-Linear Integral-PDEs without the Assumption {\bf (BD)}}

\bt\label{PDE existence}
Suppose that $\beta\in(1,2)$ and the assumptions {\bf (L)}, {\bf (FD)}, {\bf (FE)} and {\bf (BL)} are satisfied. Then
\begin{enumerate}
\item
There exists a unique mild solution $v$ of Eq. (\ref{PDE}) given by the formula
\be\label{FK}
v(t,x)=Y(t,t,x).
\ee
\item
Moreover, there exists a constant $C>0$ such that
\be
|v(t,x)|&\leq& C(1+|x|^{\mu}),\label{estimate PDE}\\
|\nabla_xv(t,x)|&\leq& C(T-t)^{-\frac{1}{\beta}}(1+|x|^{\mu}),\quad x\in\mR^d,~t\in[0,T].\label{gradient estimate PDE}
\ee
\item
Furthermore, for any $x,h\in\mR^d$, $t\in[0,T]$, we have
\be\label{derivative PDE}
\begin{split}
&\nabla_xv(t,x)h\\
=&-\mE[\int_t^T\psi(s,X(s,t,x),v(s,X(s,t,x)),\int_{\cO}(v(s,X(s,t,x)+\sigma(s,X(s,t,x))u)\\
&-v(s,X(s,t,x)))l(u)\nu(du))\cdot[\delta_{\sharp}(F_{(s-t)^{\frac{1}{\beta}}}^*(s)\cdot h)\cdot G_{(s-t)^{\frac{1}{\beta}}}^{-1}(s)\\
&-\widehat{\mE}[F_{(s-t)^{\frac{1}{\beta}}}^*(s)\cdot h\cdot (G_{(s-t)^{\frac{1}{\beta}}}^{-1}(s))^{\sharp}]]ds]\\
&+\mE[\phi(X(T,t,x))\cdot[\delta_{\sharp}(F_{(T-t)^{\frac{1}{\beta}}}^*(T)\cdot h)\cdot G_{(T-t)^{\frac{1}{\beta}}}^{-1}(T)\\
&-\widehat{\mE}[F_{(T-t)^{\frac{1}{\beta}}}^*(T)\cdot h\cdot (G_{(T-t)^{\frac{1}{\beta}}}^{-1}(T))^{\sharp}]]].
\end{split}
\ee
\end{enumerate}
\et

\br
It is worth noting that an existence and uniqueness result of a mild solution for a non-linear path-dependent partial integro-differential equation was established in \cite{CPO}, even without assuming any nondegeneracy on $\sigma$ or requiring any smoothing property of $P_{t,\tau}$. The main contribution of our article is to prove that, under the assumption of nondegeneracy on $\sigma$, the existence and uniqueness of a solution for nonlocal quasi-linear integral partial differential equations can be obtained. These solutions are differentiable with respect to the space variable, even if the initial datum and coefficients of the equation are not differentiable.
\er

Now we are ready to prove Theorem \ref{PDE existence}. Since the function $\phi$ and $\psi$ do not satisfied the assumption {\bf (BD)}, 
Theorem \ref{gradient estimate} can not be applied directly. Therefore the following lemma will play an important role in the proof of Theorem \ref{PDE existence}. For its detailed proof, the readers is referred to \cite[Lemma 4.3]{FT}.

\bl\label{mollifier}
Let $\phi$ and $\psi$ satisfy the assumption {\bf (BD)}, then for any $n\in\mN^*$, there exist $\phi_n:\mR^d\rightarrow\mR$ and $\psi_n:[0,T]\times\mR^d\times\mR\times\mR\rightarrow\mR$ such that
\begin{enumerate}
\item
$\phi_n$, $\psi_n$ satisfy the assumption {\bf (BD)} with constants $C_{BDn}$ and $m_n$ depending on $n$;
\item
$\phi_n$, $\psi_n$ satisfy the assumption {\bf (BL)} with constants $C_{BL}$ and $\mu$ independent of $n$;
\item
$\phi_n\rightarrow\phi$, $\psi_n\rightarrow\psi$ pointwise as $n\rightarrow\infty$.
\end{enumerate}
\el

\begin{proof}[Completion of the proof of Theorem \ref{PDE existence}]

The uniqueness statement in Theorem \ref{PDE existence} can be proved exactly as in Theorem \ref{PDE existence BD}. To complete the proof of Theorem \ref{PDE existence}, it remains to show that the function $v$ in Eq. (\ref{FK}) is the required solution, and that the estimates (\ref{estimate PDE}), (\ref{gradient estimate PDE}) and the formula (\ref{derivative PDE}) hold.

{\bf Step 1:}

It suffices to prove the result for $t$ varying in an interval $[T-\delta,T]$ with $\delta$ sufficiently small. 
We notice that $v(t,x)$ is deterministic and Lemma \ref{flow} yields, $\mP$-a.s.,
\ce
Y(\tau,t,x)=v(\tau,X(\tau,t,x)),\quad\tau\in[t,T],
\de
where $X$ is the solution to Eq. (\ref{SDE}).
Now let $\psi_n$, $\phi_n$ be defined in Lemma \ref{mollifier}, and hence $\{Y^n(\tau,t,x),Z^n(\tau,t,x,u);\tau\in[t,T],x\in\mR^d,u\in\cO\}$, satisfies the following backward SDEs with jumps:
\ce
&&Y^n(\tau,t,x)+\int_{\tau}^T\int_{\cO}Z^n(s-,t,x,u)\widetilde{N}(ds,du)\\
&=&-\int_{\tau}^T\psi_n(s,X(s,t,x),Y^n(s,t,x),\int_{\cO}Z^n(s,t,x,u)l(u)\nu(du))ds+\phi_n(X(T,t,x)).
\de
We also set $v_n(t,x)=Y^n(t,t,x)$. Then it is obvious that $v_n(t,x)$ is deterministic, and it can be seen from (\ref{BSDE estimate}) that
\ce
|v_n(t,x)|\leq C(1+|x|^{\mu}),\quad x\in\mR^d,\quad t\in[0,T],
\de
where $C$ is independent of $n$. Using Theorem \ref{PDE existence BD}, it follows that $v_n$ belong to $\cC^{0,1}([0,T]\times\mR^d;\mR)$, and $\mP$-a.s., 
\ce
Y^n(\tau,t,x)&=&v_n(\tau,X(\tau,t,x)),\quad\forall\tau\in[t,T],\\
Z^n(\tau,t,x,u)&=&v_n(\tau,X(\tau,t,x)+\sigma(\tau,X(\tau,t,x))u)-v_n(\tau,X(\tau,t,x)),\quad\forall\text{a.e.}~d\tau\nu(du).
\de
Therefore, it holds that
\be\label{PDE approximation}
\begin{split}
&v_n(t,x)\\
=&-\mE[\int_t^T\psi_n(\tau,X(\tau,t,x),v_n(\tau,X(\tau,t,x)),\int_{\cO}Z(\tau,x,t,u)l(u)\nu(du))d\tau]\\
&+\mE[\phi_n(X(\tau,t,x))]\\
=&-\mE[\int_t^T\psi_n(\tau,X(\tau,t,x),v_n(\tau,X(\tau,t,x)),\\
&\int_{\cO}(v_n(\tau,X(\tau,t,x)+\sigma(\tau,X(\tau,t,x))u)-v_n(\tau,X(\tau,t,x)))l(u)\nu(du))d\tau]\\
&+\mE[\phi_n(X(\tau,t,x))]\\
=&-\int_t^TP_{t,\tau}[\psi_n(\tau,\cdot,v_n(\tau,\cdot),\int_{\cO}(v_n(\tau,\cdot+\sigma(\tau,\cdot)u)-v_n(\tau,\cdot))l(u)\nu(du))](x)d\tau\\
&+P_{t,T}[\phi_n](x).
\end{split}
\ee

{\bf Step 2:}

For simplicity, in the rest of the proof, we denote
\ce
(X_{\tau},Y_{\tau},Z_{\tau}(u))=(X(\tau,t,x),Y(\tau,t,x),Z(\tau,t,x,u))
\de
and
\ce
(X_{\tau}^n,Y_{\tau}^n,Z_{\tau}^n(u))=(X^n(\tau,t,x),Y^n(\tau,t,x),Z^n(\tau,t,x,u)).
\de
It is easy to see that
\ce
&&Y_{\tau}-Y_{\tau}^n+\int_{\tau}^T\int_{\cO}Z_{s-}(u)\widetilde{N}(ds,du)-\int_{\tau}^T\int_{\cO}Z_{s-}^n(u)\widetilde{N}(ds,du)\\
&=&-\int_{\tau}^T(\psi(s,X_s,Y_s,\int_{\cO}Z_s(u)l(u)\nu(du))-\psi_n(s,X_s,Y_s,\int_{\cO}Z_s(u)l(u)\nu(du)))ds\\
&&-\int_{\tau}^T(\psi_n(s,X_s,Y_s,\int_{\cO}Z_s(u)l(u)\nu(du))-\psi_n(s,X_s^n,Y_s^n,\int_{\cO}Z_s^n(u)l(u)\nu(du)))ds\\
&&+\phi(X_T)-\phi_n(X_T).
\de
This yields by Theorem \ref{BSDE existence} that
\ce
&&\mE[\sup_{\tau\in[t,T]}|Y_{\tau}-Y_{\tau}^n|]+\mE[\int_t^T\int_{\cO}|Z_{\tau}-Z_{\tau}^n|^2\nu(du)d\tau]\\
&\leq&C\mE[|\phi(X_T)-\phi_n(X_T)|^2]\\
&&+C(T-t)\mE[\int_t^T|\psi(s,X_s,Y_s,\int_{\cO}Z_s(u)l(u)\nu(du))-\psi_n(s,X_s,Y_s,\int_{\cO}Z_s(u)l(u)\nu(du))|^2]\\
&&+C(T-t)\mE[\int_t^T[|Y_s-Y_s^n|^2+|\int_{\cO}(Z_s(u)-Z_s^n(u))l(u)\nu(du)|^2]ds]\\
&\leq&C\mE[|\phi(X_T)-\phi_n(X_T)|^2]\\
&&+C(T-t)\mE[\int_t^T|\psi(s,X_s,Y_s,\int_{\cO}Z_s(u)l(u)\nu(du))-\psi_n(s,X_s,Y_s,\int_{\cO}Z_s(u)l(u)\nu(du))|^2]\\
&&+C(T-t)\mE[\int_t^T[|Y_s-Y_s^n|^2+\int_{\cO}|Z_s(u)-Z_s^n(u)|^2\nu(du)]ds]\\
&=:&J_1^n+J_2^n+J_3^n,
\de
where we have used the following fact:
\ce
&&|\int_{\cO}(Z_s(u)-Z_s^n(u))l(u)\nu(du)|^2\\
&=&(\int_{\cO}l^2(u)\nu(du))^2(\int_{\cO}|\frac{Z_s(u)-Z_s^n(u)}{l(u)}|\frac{l^2(u)\nu(du)}{\int_{\cO}l^2(u)\nu(du)})^2\\
&\leq&\int_{\cO}l^2(u)\nu(du)\int_{\cO}|Z_s(u)-Z_s^n(u)|^2\nu(du).
\de
Using Lebesgue dominate convergence theorem, it is easy to see that
\ce
J_1^n,J_2^n\rightarrow0,\quad\text{as}\quad n\rightarrow\infty.
\de
Then using the fact that $T-t\leq\delta$ small enough, we have
\ce
\mE[\sup_{\tau\in[t,T]}|Y_{\tau}-Y_{\tau}^n|]+\mE[\int_t^T\int_{\cO}|Z_{\tau}(u)-Z_{\tau}^n(u)|^2\nu(du)d\tau]\rightarrow0,~\text{as}~n\rightarrow\infty.
\de
This implies that
\be\label{BSDE estimate 2}
\begin{split}
&\mE[\sup_{\tau\in[t,T]}|Y_{\tau}-Y_{\tau}^n|]\\
&+\mE[\int_t^T\int_{\cO}|Z_{\tau}(u)-(v_n(\tau,X_{\tau}+\sigma(\tau,X_{\tau})u)-v_n(\tau,X_{\tau}))|^2\nu(du)d\tau]\\
&\rightarrow0,~\text{as}~n\rightarrow\infty.
\end{split}
\ee
In particular, it is obvious that $v_n(t,x)\rightarrow v(t,x)$ for every $t\in[0,T]$, $x\in\mR^d$.

{\bf Step 3:}

Applying the estimate in Theorem \ref{gradient estimate}, we obtain
\be\label{gradient estimate approximation}
|\nabla_xv_n(t,x)|\leq C(T-t)^{-1/\beta}(1+|x|)^{\mu},\quad\forall t\in[0,T],~\forall x\in\mR^d.
\ee
This implies by the mean value theorem, the assumption \textbf{(FE)} and Theorem \ref{FSDE existence} that $\mP$-a.s. for a.e. $s\in[t,T]$,
\be\label{BSDE estimate 3}
\begin{split}
&\int_{\cO}|Z_s^n(u)||l(u)|\nu(du)\\
=&\int_{\cO}|v_n(s,X(s,t,x)+\sigma(s,X(s,t,x))u)-v_n(s,X(s,t,x))||l(u)|\nu(du)\\
=&\int_{\cO}|\nabla_x v_n(s,X(s,t,x)+(1-\theta)\sigma(s,X(s,t,x))u)\sigma(s,X(s,t,x))||ul(u)|\nu(du)\\
\leq&C_{FB}\int_{\cO}|\nabla_x v_n(s,X(s,t,x)+(1-\theta)\sigma(s,X(s,t,x))u)||ul(u)|\nu(du)\\
\leq&C(T-s)^{-1/\beta}(1+|X(s,t,x)|+(1-\theta)|\sigma(s,X(s,t,x))|)^{\mu}\\
\leq&C(T-s)^{-1/\beta}(1+|X(s,t,x)|)^{\mu},
\end{split}
\ee
where $0<\theta<1$, and hence by (\ref{BSDE estimate 2})
\be\label{BSDE estimate 4}
\int_{\cO}|Z_s(u)||l(u)|\nu(du)\leq C(T-s)^{-1/\beta}(1+|X(s,t,x)|)^{\mu}.
\ee

{\bf Step 4:}

Now we define
\ce
&&V(t,x)h\\
&:=&-\mE[\int_t^T\psi(s,X(s,t,x),Y(s,t,x),\int_{\cO}Z(s,t,x,u)l(u)\nu(du))\cdot[\delta_{\sharp}(F_{(s-t)^{\frac{1}{\beta}}}^*(s)\cdot h)\cdot G_{(s-t)^{\frac{1}{\beta}}}^{-1}(s)\\
&&-\widehat{\mE}[F_{(s-t)^{\frac{1}{\beta}}}^*(s)\cdot h\cdot (G_{(s-t)^{\frac{1}{\beta}}}^{-1}(s))^{\sharp}]]ds]\\
&&+\mE[\phi(X(T,t,x))\cdot[\delta_{\sharp}(F_{(T-t)^{\frac{1}{\beta}}}^*(T)\cdot h)\cdot G_{(T-t)^{\frac{1}{\beta}}}^{-1}(T)\\
&&-\widehat{\mE}[F_{(T-t)^{\frac{1}{\beta}}}^*(T)\cdot h\cdot (G_{(T-t)^{\frac{1}{\beta}}}^{-1}(T))^{\sharp}]]],\quad t\in[0,T],~x\in\mR^d.
\de
The following result can be proved:
\ce
|V(t,x)-\nabla_xv_n(t,x)|\rightarrow0,\quad t\in[0,T],~x\in\mR^d,~\text{as}~n\rightarrow\infty.
\de
Indeed, it can be see from Theorem \ref{Bismut formula} that
\ce
&&\nabla v_n(t,x)h\\
&=&-\mE[\int_t^T\psi_n(s,X(s,t,x),Y(s,t,x),\int_{\cO}Z(s,t,x,u)l(u)\nu(du))\\
&&\cdot[\delta_{\sharp}(F_{(s-t)^{\frac{1}{\beta}}}^*(s)\cdot h)\cdot G_{(s-t)^{\frac{1}{\beta}}}^{-1}(s)-\widehat{\mE}[F_{(s-t)^{\frac{1}{\beta}}}^*(s)\cdot h\cdot (G_{(s-t)^{\frac{1}{\beta}}}^{-1}(s))^{\sharp}]]ds]\\
&&+\mE[\phi_n(X(T,t,x))\cdot[\delta_{\sharp}(F_{(T-t)^{\frac{1}{\beta}}}^*(T)\cdot h)\cdot G_{(T-t)^{\frac{1}{\beta}}}^{-1}(T)\\
&&-\widehat{\mE}[F_{(T-t)^{\frac{1}{\beta}}}^*(T)\cdot h\cdot (G_{(T-t)^{\frac{1}{\beta}}}^{-1}(T))^{\sharp}]]].
\de
Recall that
\ce
U_s^h=\delta_{\sharp}(F_{(s-t)^{\frac{1}{\beta}}}^*(s)\cdot h)\cdot G_{(s-t)^{\frac{1}{\beta}}}^{-1}(s)-\widehat{\mE}[F_{(s-t)^{\frac{1}{\beta}}}^*(s)\cdot h\cdot (G_{(s-t)^{\frac{1}{\beta}}}^{-1}(s))^{\sharp}].
\de
Then we have
\ce
&&V(t,x)h-\nabla v_n(t,x)h\\
&=&-\mE[\int_t^T(\psi(s,X(s,t,x),Y(s,t,x),\int_{\cO}Z(s,t,x,u)l(u)\nu(du))\cdot U_s^h\\
&&-\psi_n(s,X(s,t,x),Y^n(s,t,x),\int_{\cO}Z^n(s,t,x,u)l(u)\nu(du))\cdot U_s^h)ds]\\
&&+\mE[(\phi(X(T,t,x))-\phi_n(X(T,t,x)))\cdot U_T^h]\\
&=&-\mE[\int_t^T(\psi(s,X(s,t,x),Y(s,t,x),\int_{\cO}Z(s,t,x,u)l(u)\nu(du))\cdot U_s^h\\
&&-\psi_n(s,X(s,t,x),Y(s,t,x),\int_{\cO}Z(s,t,x,u)l(u)\nu(du))\cdot U_s^h)ds]\\
&&-\mE[\int_t^T(\psi_n(s,X(s,t,x),Y(s,t,x),\int_{\cO}Z(s,t,x,u)l(u)\nu(du))\cdot U_s^h\\
&&-\psi_n(s,X(s,t,x),Y^n(s,t,x),\int_{\cO}Z^n(s,t,x,u)l(u)\nu(du))\cdot U_s^h)ds]\\
&&+\mE[(\phi(X(T,t,x))-\phi_n(X(T,t,x)))\cdot U_T^h].
\de
This implies by Theorem \ref{BSDE existence}, H\"{o}lder's inequality and Lemma \ref{Bismut estimate} that
\ce
&&|V(t,x)-\nabla v_n(t,x)|\\
&\leq&C\int_t^T\mE[|\psi(s,X(s,t,x),v(s,X(s,t,x)),\int_{\cO}Z(s,t,x,u)l(u)\nu(du))\\
&&-\psi_n(s,X(s,t,x),v(s,X(s,t,x)),\int_{\cO}Z(s,t,x,u)l(u)\nu(du))|^2]^{1/2}(s-t)^{-1/\beta}ds\\
&&+C\int_t^T\mE[|\psi_n(s,X(s,t,x),v(s,X(s,t,x)),\int_{\cO}Z(s,t,x,u)l(u)\nu(du))\\
&&-\psi_n(s,X(s,t,x),v_n(s,X(s,t,x)),\int_{\cO}Z^n(s,t,x,u)l(u)\nu(du))|^2]^{1/2}(s-t)^{-1/\beta}ds\\
&&+C\int_t^T\mE[|\phi(X(T,t,x))-\phi_n(X(T,t,x))|^2]^{1/2}(T-t)^{-1/\beta}ds\\
&\leq&C\int_t^T\mE[|\psi(s,X(s,t,x),v(s,X(s,t,x)),\int_{\cO}Z(s,t,x,u)l(u)\nu(du))\\
&&-\psi_n(s,X(s,t,x),v(s,X(s,t,x)),\int_{\cO}Z(s,t,x,u)l(u)\nu(du))|^2]^{1/2}(s-t)^{-1/\beta}ds\\
&&+C\int_t^T\mE[|v(s,X(s,t,x))-v_n(s,X(s,t,x))|^2]^{1/2}(s-t)^{-1/\beta}ds\\
&&+C\int_t^T\mE[|\int_{\cO}(Z(s,t,x,u)-(v_n(s,X(s,t,x)+\sigma(s,X(s,t,x))u)\\
&&-v_n(s,X(s,t,x))))l(u)\nu(du)|^2]^{1/2}(s-t)^{-1/\beta}ds\\
&&+C\int_t^T\mE[|\phi(X(T,t,x))-\phi_n(X(T,t,x))|^2]^{1/2}(T-t)^{-1/\beta}ds\\
&=:&C(K_1^n+K_2^n+K_3^n+K_4^n).
\de
It can be seen from Lebesgue dominated convergence theorem that
\ce
K_1^n,K_2^n,K_4^n\rightarrow0,\quad\text{as}\quad n\rightarrow\infty.
\de
Now we prove that
\ce
K_3^n\rightarrow0,\quad\text{as}\quad n\rightarrow\infty.
\de
Suppose the contrary, i.e., for some $\epsilon>0$ and some subsequence $n_k$,
\ce
&&K_3^{n_k}\\
&=&\int_t^T\mE[|\int_{\cO}(Z(s,t,x,u)-(v_{n_k}(s,X(s,t,x)+\sigma(s,X(s,t,x))u)\\
&&-v_{n_k}(s,X(s,t,x))))l(u)\nu(du)|^2]^{1/2}(s-t)^{-1/\beta}ds\\
&\geq&\epsilon.
\de
However, using (\ref{BSDE estimate 2}), we can extract a subsequence $(n_{k_j})$ such that
\ce
&&\mE[|\int_{\cO}(Z(s,t,x,u)-(v_{n_{k_j}}(s,X(s,t,x)+\sigma(s,X(s,t,x))u)\\
&&-v_{n_{k_j}}(s,X(s,t,x))))l(u)\nu(du)|^2]\rightarrow0,\quad\text{as}\quad n\rightarrow\infty,
\de
and hence by (\ref{BSDE estimate 3})-(\ref{BSDE estimate 4}), Theorem \ref{FSDE existence} and Lebesgue dominated convergence theorem, a contradiction can be derived. Then the desired result is obtained.

{\bf Step 5:}

It is implied by the continuous dependence on $(t,x)$ of $X$ that $(t,x)\mapsto V(t,x)$ is also continuous. Therefore $v(t,\cdot)$ admits directional derivatives in every direction $h\in\mR^d$ at every point $x\in\mR^d$, equal to $V(t,x)h$. This implies that $v\in\cC^{0,1}([0,T]\times\mR^d;\mR)$ and $\nabla_xv=V$. Using the fact that for a.e. $(\tau,u)$, $\mP$-a.s.
\ce
Z_{\tau}^n(u)=v_n(\tau,X_{\tau}+\sigma(\tau,X_{\tau})u)-v_n(\tau,X_{\tau}),
\de
we have
\ce
Z_{\tau}(u)=v(\tau,X_{\tau}+\sigma(\tau,X_{\tau})u)-v(\tau,X_{\tau}).
\de
Then the formula (\ref{derivative PDE}) can be obtained from the definition of $V$. The fact that $v$ solves Eq. (\ref{PDE}) come from (\ref{PDE approximation}) by letting $n\rightarrow\infty$.

The estimate (\ref{estimate PDE}) is a direct result of (\ref{BSDE estimate}). Using (\ref{gradient estimate approximation}), it is also easy to see that the estimate (\ref{gradient estimate PDE}) is satisfied.
Then the proof is complete.
\end{proof}

\begin{appendix}
\section{Proof of Theorem \ref{BSDE existence}}\label{appA}
The proof of the first three conclusion of Theorem \ref{BSDE existence} is standard, and the readers are refered to e.g. \cite{LW1,LW2,L,D}. Therefore it is sufficient to prove the fourth conclusion. We denote by
\begin{enumerate}
\item
$\cH_{\mD,\cP}$ the subvector space of predictable processes in $\cH_{\mD}$;
\item
$\cH_{\mD,\nu,\cP}$ the subvector space of predictable processes in $\cH_{\mD,\nu}$;
\item
$\cH_{\mD\otimes\mathbf{d},\cP}$ the set of real valued processes $H$ defined on $[0,T]\times\Omega\times\cO$ which are predictable and belong to $L^2([0,T];\mD\times\mathbf{d})$, i.e., such that
\ce
\|H\|^2_{\cH_{\mD\otimes\mathbf{d},\cP}}&=&\mE[\int_0^T\int_{\cO}|H(t,u)|^2\nu(du)dt]\\
&&+\int_0^T\int_{\cO}\cE(H(t,u))\nu(du)dt+\mE[\int_0^Te(H(t,\cdot))dt]<\infty.
\de
\end{enumerate}
We also denote by $\cH_{\mD,\cP}^d$, $\cH_{\mD,\nu,\cP}^d$ and $\cH_{\mD\otimes\mathbf{d},\cP}^d$ the space of $\mR^d$-valued processes such that each coordinate belongs, respectively, to $\cH_{\mD,\cP}$, $\cH_{\mD,\nu,\cP}$ and $\cH_{\mD\otimes\mathbf{d},\cP}$ and we equip it with the standard norm of product topology.

We shall need the following results which can be obtained using the same argument as in the proof \cite[Proposition 8.2]{BD2}.

\bl\label{gradient of stochastic integration}
Let $H\in\cH_{\mD\otimes\mathbf{d},\cP}$, $G\in\cH_{\mD,\cP}$ and $F\in\cH_{\mD,\nu,\cP}$, then
\begin{enumerate}
\item
The process
\ce
X(t)=\int_0^t\int_{\cO}H(s,\omega,u)\widetilde{N}(ds,du),\quad t\in[0,T]
\de
is a square integrable martingale which belongs to $\cH_{\mD}$, and
\ce
X^{\sharp}(t,\omega,\hat{\omega})&=&\int_0^t\int_{\cO}H^{\sharp}(s,\omega,u,\hat{\omega})\widetilde{N}(ds,du)\\
&&+\int_0^t\int_{\cO}\int_RH^{\flat}(s,\omega,u,r)N\odot\rho(ds,du,dr),\quad\forall t\in[0,T].
\de
Moreover,
\ce
\|X(t)\|_{\mD}&\leq&\sqrt{2}\|H\|_{\cH_{\mD\otimes\mathbf{d},\cP}},\quad\forall t\in[0,T],\\
\|X\|_{\cH_{\mD}}&\leq&\sqrt{2T}\|H\|_{\cH_{\mD\otimes\mathbf{d},\cP}}.
\de
\item
The process
\ce
Y(t)=\int_0^tG(s,\omega)ds,\quad\forall t\in[0,T]
\de
is a square integrable semimartingale which belongs to $\cH_{\mD}$, and
\ce
Y^{\sharp}(t,\omega,\hat{\omega})=\int_0^t G^{\sharp}(s,\omega,\hat{\omega})ds,\quad\forall t\in[0,T].
\de
Moreover,
\ce
\|Y(t)\|_{\mD}&\leq&C\|G\|_{\cH_{\mD,\cP}},\quad\forall t\in[0,T],\\
\|Y\|_{\cH_{\mD}}&\leq&C\sqrt{T}\|G\|_{\cH_{\mD,\cP}}.
\de
\item
The process
\ce
Z(t)=\int_0^t\int_{\cO}F(s,u,\omega)\nu(du)ds,\quad\forall t\in[0,T]
\de
is a square integrable semimartingale which belongs to $\cH_{\mD}$, and
\ce
Z^{\sharp}(t,\omega,\hat{\omega})=\int_0^t\int_{\cO} F^{\sharp}(s,u,\omega,\hat{\omega})\nu(du)ds,\quad\forall t\in[0,T].
\de
Moreover,
\ce
\|Z(t)\|_{\mD}&\leq&C\|F\|_{\cH_{\mD,\nu,\cP}},\quad\forall t\in[0,T],\\
\|Z\|_{\cH_{\mD}}&\leq&C\sqrt{T}\|F\|_{\cH_{\mD,\nu,\cP}}.
\de
\end{enumerate}
\el

The following commutation relation between the Malliavin derivative of a random variable and its conditional expectation with respect to $\mathcal{F}_s$ will  also be used.

\bl\label{conditional}
Suppose that $\{g(t,u,r);t\in[0,T],u\in\cO,r\in R\}$ satisfies
\ce
\int_0^T\int_{\cO\times R}g^2(t,u,r)N(dt,du)\rho(r)<\infty.
\de
Let $F\in\mD$ and $s\in[0,T]$. Then $\mE[F|\mathcal{F}_s]\in\mD$
and
\ce
\widehat{\mE}[(\mE[F|\mathcal{F}_s])^{\sharp}G(T)]=\widehat{\mE}[\mE[F^{\sharp}|\mathcal{F}_s]G(s)],
\de
where $G(s)=\int_0^T\int_{\cO\times R}g(s',u,r)1_{[0,s]}(s')N\odot\rho(ds',du,dr)$.

\begin{proof}
Let $\cE$ be the linear span of all real parts of functions
\ce
\{e^{-\int_0^T\int_{\cO}f(s',u)\widetilde{N}(ds',du)};f\geq0,f\in\mathbf{d}\cap L^1\cap L^{\infty},\int_0^T\int_{\cO}f^2(s',u)\nu(du)ds'<\infty\}.
\de
It can be seen from Lemma \cite[Lemma 4.39]{BD3} that $\cE$ is total in $\mD$.

Let us first prove Lemma \ref{conditional} when $F=e^{-\int_0^T\int_{\cO}f(s',u)\widetilde{N}(ds',du)}\in\cE$. Using It\^{o}'s formula, it follows that
\ce
&&e^{-\int_0^t\int_{\cO}f(s',u)\widetilde{N}(ds',du)}\\
&=&e^{-\int_0^s\int_{\cO}f(s',u)\widetilde{N}(ds',du)}\\
&&+\int_s^t\int_{\cO}(e^{-\int_0^{s'-}\int_{\cO}f(s'',u)\widetilde{N}(ds'',du)-f(s',u)}-e^{-\int_0^{s'-}\int_{\cO}f(s'',u)\widetilde{N}(ds'',du)})\widetilde{N}(ds',du)\\
&&+\int_s^t\int_{\cO}(e^{-\int_0^{s'-}\int_{\cO}f(s'',u)\widetilde{N}(ds'',du)-f(s',u)}-e^{-\int_0^{s'-}\int_{\cO}f(s'',u)\widetilde{N}(ds'',du)}\\
&&+f(s',u)e^{-\int_0^{s'-}\int_{\cO}f(s'',u)\widetilde{N}(ds'',du)})\nu(du)ds'.
\de
Taking conditional expectation with respect to $\mathcal{F}_s$ of both sides of this identity and setting $C(t)=\mE[e^{-\int_0^t\int_{\cO}f(s',u)\widetilde{N}(ds',du)}|\mathcal{F}_s]$ yields
\ce
C(t)&=&C(s)+\int_s^t\int_{\cO}(e^{-f(s',u)}C(s'-)-C(s'-)+f(s',u)C(s'-))\nu(du)ds'.
\de
This implies that
\ce
C(t)=C(s)e^{\int_s^t\int_{\cO}(e^{-f(s',u)}-1+f(s',u))\nu(du)ds'}.
\de
Then
\ce
\mE[e^{-\int_0^T\int_{\cO}f(s',u)\widetilde{N}(ds',du)}|\mathcal{F}_s]=e^{-\int_0^s\int_{\cO}f(s',u)\widetilde{N}(ds',du)}e^{-\int_s^T\int_{\cO}(e^{f(s',u)}-1+f(s',u))\nu(du)ds'}\in\mD,
\de
and hence by the chain rule,
\ce
&&\widehat{\mE}[(\mE[e^{-\int_0^T\int_{\cO}f(s',u)\widetilde{N}(ds',du)}|\mathcal{F}_s])^{\sharp}G(T)]\\
&=&\widehat{\mE}[-\int_0^s\int_{\cO\times\mR}f^{\flat}(s',u,r)N\odot\rho(ds',du,dr)e^{-\int_0^s\int_{\cO}f(s',u)\widetilde{N}(ds',du)}e^{\int_s^t\int_{\cO}(e^{-f(s',u)}-1+f(s',u))\nu(du)ds'}\\
&&\int_0^T\int_{\cO\times R}g(s',u,r)N\odot\rho(ds',du,dr)]\\
&=&-\int_0^s\int_{\cO\times R}f^{\flat}(s',u,r)g(s',u,r)N(ds',du)\rho(dr)\\
&&e^{-\int_0^s\int_{\cO}f(s',u)\widetilde{N}(ds',du)}e^{-\int_s^t\int_{\cO}(e^{f(s',u)}-1+f(s',u))\nu(du)ds'}.
\de
On the other hand, it is obvious that
\ce
(e^{-\int_0^T\int_{\cO}f(s',u)\widetilde{N}(ds',du)})^{\sharp}=-\int_0^T\int_{\cO\times\mR}f^{\flat}(s',u,r)N\odot\rho(ds',du,dr)e^{-\int_0^T\int_{\cO}f(s',u)\widetilde{N}(ds',du)}.
\de
Then using It\^{o}'s product formula, we have
\be\label{product}
\begin{split}
&-\int_0^t\int_{\cO\times\mR}f^{\flat}(s',u,r)N\odot\rho(ds',du,dr)e^{-\int_0^t\int_{\cO}f(s',u)\widetilde{N}(ds',du)}\\
=&-\int_0^s\int_{\cO\times\mR}f^{\flat}(s',u,r)N\odot\rho(ds',du,dr)e^{-\int_0^s\int_{\cO}f(s',u)\widetilde{N}(ds',du)}\\
&-\int_s^t\int_0^{s'-}\int_{\cO\times\mR}f^{\flat}(s'',u,r)N\odot\rho(ds'',du,dr)d(e^{-\int_0^{s'}\int_{\cO}f(s'',u)\widetilde{N}(ds'',du)})\\
&-\int_s^te^{-\int_0^{s'-}\int_{\cO}f(s'',u)\widetilde{N}(ds'',du)}d(\int_0^{s'}\int_{\cO\times\mR}f^{\flat}(s'',u,r)N\odot\rho(ds'',du,dr))\\
=&-\int_0^s\int_{\cO\times\mR}f^{\flat}(s',u,r)N\odot\rho(ds',du,dr)e^{-\int_0^s\int_{\cO}f(s',u)\widetilde{N}(ds',du)}\\
&-\int_s^t\int_0^{s'-}\int_{\cO\times\mR}f^{\flat}(s'',u,r)N\odot\rho(ds'',du,dr)\\
&d(\int_s^{s'}\int_{\cO}(e^{-\int_0^{s''-}\int_{\cO}f(s''',u)\widetilde{N}(ds''',du)-f(s'',u)}\\
&-e^{-\int_0^{s''-}\int_{\cO}f(s''',u)\widetilde{N}(ds''',du)})\widetilde{N}(ds'',du)\\
&+\int_s^{s'}\int_{\cO}(e^{-\int_0^{s''-}\int_{\cO}f(s''',u)\widetilde{N}(ds''',du)-f(s'',u)}-e^{-\int_0^{s''-}\int_{\cO}f(s''',u)\widetilde{N}(ds''',du)}\\
&+f(s'',u)e^{-\int_0^{s''-}\int_{\cO}f(s''',u)\widetilde{N}(ds''',du)})\nu(du)ds'')\\
&-\int_s^te^{-\int_0^{s'-}\int_{\cO}f(s'',u)\widetilde{N}(ds'',du)}d(\int_0^{s'}\int_{\cO\times\mR}f^{\flat}(s'',u,r)N\odot\rho(ds'',du,dr)).
\end{split}
\ee
Now setting
\ce
D(t)=-\widehat{\mE}[\mE[\int_0^t\int_{\cO\times\mR}f^{\flat}(s',u,r)N\odot\rho(ds',du,dr)e^{-\int_0^t\int_{\cO}f(s',u)\widetilde{N}(ds',du)}|\mathcal{F}_s]G(s)],
\de
and hence it follows from (\ref{product}) that
\ce
D(t)=-\int_0^s\int_{\cO\times\mR}f^{\flat}(s',u,r)\widetilde{N}(ds',du)\rho(dr)+\int_s^t\int_{\cO}D(s'-)(e^{-f(s',u)}-1+f(s',u))\nu(du)ds'.
\de
Therefore
\ce
D(t)=-\int_0^s\int_{\cO\times\mR}f^{\flat}(s',u,r)\widetilde{N}(ds',du)\rho(dr)e^{\int_s^t\int_{\cO}(e^{-f(s',u)}-1+f(s',u))\nu(du)ds'}.
\de
This yields that
\ce
&&\widehat{\mE}[\mE[(e^{-\int_0^T\int_{\cO}f(s',u)\widetilde{N}(ds',du)})^{\sharp}|\mathcal{F}_s]G_s]\\
&=&-\int_0^s\int_{\cO\times\mR}f^{\flat}(s',u,r)\widetilde{N}(ds',du)\rho(dr)e^{-\int_0^s\int_{\cO}f(s',u)\widetilde{N}(ds',du)}e^{\int_s^t\int_{\cO}(e^{-f(s',u)}-1+f(s',u))\nu(du)ds'}.
\de
Then Lemma \ref{conditional} is proved when $F\in\cE$.

Let now $F\in\mD$ and $(F_n,n\in\mN^*)\subset\cE$ be such that
\ce
\lim_{n\rightarrow\infty}F_n=F\quad&\text{in}&\quad L^2(\mP),\\
\lim_{n\rightarrow\infty}F_n^{\sharp}=F^{\sharp}\quad&\text{in}&\quad L^2(\mP\times\widehat{\mP}).
\de
Using the continuity of the conditional expectation, we have
\ce
\lim_{n\rightarrow\infty}\mE[F_n|\mathcal{F}_s]=\mE[F|\mathcal{F}_s]\quad&\text{in}&\quad L^2(\mP),\\
\lim_{n\rightarrow\infty}\mE[F_n^{\sharp}|\mathcal{F}_s]=\mE[F^{\sharp}|\mathcal{F}_s]\quad&\text{in}&\quad L^2(\mP\times\widehat{\mP}).
\de
Therefore applying
\ce
\widehat{\mE}[(E[F_n|\mathcal{F}_s])^{\sharp}G(T)]=\widehat{\mE}[\mE[F_n^{\sharp}|\mathcal{F}_s]G(s)],
\de
it follows that
\ce
\lim_{n\rightarrow\infty}\widehat{\mE}[(\mE[F_n|\mathcal{F}_s])^{\sharp}G(T)]=\widehat{\mE}[\mE[F^{\sharp}|\mathcal{F}_s]G(s)]\quad&\text{in}&\quad L^2(\mP).
\de
Now the desired result is obtained.
\end{proof}

\el

Using the method of Picard iteration, for every $n\in\mN$, we construct recursively a sequence $\{Y^n(\tau,t,x),Z^n(\tau,t,x,u);\tau\in[0,T],x\in\mR^d,u\in\cO\}$ of successive approximations by,
\ce
\begin{cases}
Y^{n+1}(\tau,t,x)+\int_{\tau}^T\int_{\cO}Z^{n+1}(s-,t,x,u)\widetilde{N}(ds,du)\\
=-\int_{\tau}^T\psi(s,X(s,t,x),Y^n(s,t,x),\int_{\cO}Z^n(s,t,x,u)l(u)\nu(du))ds+\phi(X(T,t,x)),\\
Y^0(\tau,t,x)=0,Z^0(\tau,t,x,u)=0.
\end{cases}
\de

The following result can be proved in the same way as the proof of \cite[Theorem 3.1.1]{D}.

\bl\label{Picard convergence}
Under the assumptions of Theorem \ref{BSDE existence}, we have
\ce
\lim_{n\rightarrow\infty}(\|Y-Y^n\|_{\cS^2}+\|Z-Z^n\|_{\cM^2})=0.
\de
\el

\bl\label{Picard Malliavin}
Under the assumptions of Theorem \ref{BSDE existence},
\ce
(Y^n(\tau,t,x),Z^n(\tau,t,x,u))\in\cH_{\mD}\times\cH_{\mD,\nu},\\
\mE\widehat{\mE}[\sup_{\tau\in[0,T]}|(Y^n(\tau,t,x))^{\sharp}|^2]+\mE\widehat{\mE}[\int_0^T\int_{\cO}|(Z^n(\tau,t,x,u))^{\sharp}|^2\nu(du)d\tau]<\infty,
\de
implies
\ce
(Y^{n+1}(\tau,t,x),Z^{n+1}(\tau,t,x,u))\in\cH_{\mD}\times\cH_{\mD,\nu},\\
\mE\widehat{\mE}[\sup_{\tau\in[0,T]}|(Y^{n+1}(\tau,t,x))^{\sharp}|^2]+\mE\widehat{\mE}[\int_0^T\int_{\cO}|(Z^{n+1}(\tau,t,x,u))^{\sharp}|^2\nu(du)d\tau]<\infty.
\de
\el

\begin{proof}

{\bf Step 1:}

By Lemma \ref{gradient of stochastic integration}, it is obvious that
\ce
\int_{\cO}Z^n(s,t,x,u)l(u)\nu(du)\in\mD,
\de
and hence by the assumption {\bf (BD)} and the chain rule, it holds that
\ce
\psi(s,X(s,t,x),Y^n(s,t,x),\int_{\cO}Z^n(s,t,x,u)l(u)\nu(du))\in\mD.
\de
Then using Lemma \ref{gradient of stochastic integration} again, for every $\tau\in[0,T]$,
\ce
\int_{\tau}^T\psi(s,X(s,t,x),Y^n(s,t,x),\int_{\cO}Z^n(s,t,x,u)l(u)\nu(du))ds\in\mD
\de
and
\ce
&&(\int_{\tau}^T\psi(s,X(s,t,x),Y^n(s,t,x),\int_{\cO}Z^n(s,t,x,u)l(u)\nu(du))ds)^{\sharp}\\
&=&\int_{\tau}^T(\nabla_x\psi(s,X(s,t,x),Y^n(s,t,x),\int_{\cO}Z^n(s,t,x,u)l(u)\nu(du))X^{\sharp}(s,t,x)\\
&&+\nabla_y\psi(s,X(s,t,x),Y^n(s,t,x),\int_{\cO}Z^n(s,t,x,u)l(u)\nu(du))(Y^n(s,t,x))^{\sharp}\\
&&+\nabla_z\psi(s,X(s,t,x),Y^n(s,t,x),\int_{\cO}Z^n(s,t,x,u)l(u)\nu(du))\int_{\cO}(Z^n(s,t,x,u))^{\sharp}l(u)\nu(du))ds.
\de

{\bf Step 2:}

By the chain rule, $\phi(X(T,t,x))\in\mD$ and
\ce
(\phi(X(T,t,x)))^{\sharp}=\nabla\phi(X(T,t,x))(X(T,t,x))^{\sharp}.
\de
In view of Lemma \ref{conditional} and
\ce
&&Y^{n+1}(\tau,t,x)\\
&=&\mE[-\int_{\tau}^T\psi(s,X(s,t,x),Y^n(s,t,x),\int_{\cO}Z^n(s,t,x,u)l(u)\nu(du))ds+\phi(X(T,t,x))|\mathcal{F}_{\tau}],
\de
for every $\tau\in[0,T]$, we have
\ce
Y^{n+1}(\tau,t,x)\in\mD.
\de
This implies that
\ce
\int_{\tau}^T\int_{\cO}Z^{n+1}(s,t,x,u)\widetilde{N}(ds,du)\in\mD.
\de

{\bf Step 3:}

For $(\alpha,y)\in[t,T]\times\mR^d$, setting
\ce
X^{(\alpha,y)}(\tau,t,x):=\varepsilon^+_{(\alpha,y)}X(\tau,t,x),\quad\nabla X^{\alpha,y}(\tau,\alpha,x):=\varepsilon^+_{(\alpha,y)}\nabla X(\tau,\alpha,x),\\
(Y^n(\tau,t,x))^{(\alpha,y)}:=\varepsilon^+_{(\alpha,y)}Y^n(\tau,t,x),\quad(\nabla Y^n(\tau,\alpha,x))^{(\alpha,y)}:=\varepsilon^+_{(\alpha,y)}\nabla Y^n(\tau,\alpha,x),\\
(Z^n(\tau,t,x,u))^{(\alpha,y)}:=\varepsilon^+_{(\alpha,y)}Z^n(\tau,t,x,u),\quad(\nabla Z^n(\tau,\alpha,x,u))^{(\alpha,y)}:=\varepsilon^+_{(\alpha,y)}\nabla Z^n(\tau,\alpha,x,u).
\de
It holds that
\ce
&&(Y^{n+1}(\tau,t,x))^{(\alpha,y)}+\int_{\tau}^T\int_{\cO}(Z^{n+1}(s-,t,x,u))^{(\alpha,y)}\widetilde{N}(ds,du)+Z^{n+1}(\alpha,t,x,y)1_{\{\tau\in[0,\alpha]\}}\\
&=&-\int_{\tau}^T\psi(s,X^{(\alpha,y)}(s,t,x),(Y ^n(s,t,x))^{(\alpha,y)},\int_{\cO}(Z^n(s,t,x,u))^{(\alpha,y)}l(u)\nu(du))ds\\
&&+\nabla\phi(X^{(\alpha,y)}(T,t,x)).
\de
Then for any $\tau\in[\alpha,T]$, we have
\ce
&&((Y^{n+1}(\tau,t,x))^{(\alpha,y)})^{\flat}+\int_{\tau}^T\int_{\cO}((Z^{n+1}(s-,t,x,u))^{(\alpha,y)})^{\flat}\widetilde{N}(ds,du)\no\\
&=&-\int_{\tau}^T(\nabla_x\psi(s,X^{(\alpha,y)}(s,t,x),(Y^n(s,t,x))^{(\alpha,y)},\no\\
&&\int_{\cO}(Z^n(s,t,x,u))^{(\alpha,y)}l(u)\nu(du))(X^{(\alpha,y)}(s,t,x))^{\flat}\no\\
&&+\nabla_y\psi(s,X^{(\alpha,y)}(s,t,x),(Y^n(s,t,x))^{(\alpha,y)},\no\\
&&\int_{\cO}(Z^n(s,t,x,u))^{(\alpha,y)}l(u)\nu(du))((Y^n(s,t,x))^{(\alpha,y)})^{\flat}\no\\
&&+\nabla_z\psi(s,X^{(\alpha,y)}(s,t,x),(Y^n(s,t,x))^{(\alpha,y)},\no\\
&&\int_{\cO}(Z^n(s,t,x,u))^{(\alpha,y)}l(u)\nu(du))\int_{\cO}((Z^n(s,t,x,u))^{(\alpha,y)})^{\flat}l(u)\nu(du))ds\no\\
&&+\nabla\phi(X^{(\alpha,y)}(T,t,x))(X^{(\alpha,y)}(T,t,x))^{\flat}.
\de
Now using (\ref{gradient}), it holds that
\ce
&&(Y^{n+1}(\tau,t,x))^{\sharp}+\int_{\tau}^T\int_{\cO}(Z^{n+1}(s-,t,x,u))^{\sharp}\widetilde{N}(ds,du)\\
&=&-\int_{\tau}^T\nabla_x\psi(s,X(s,t,x),Y^n(s,t,x),\int_{\cO}Z^n(s,t,x,u)l(u)\nu(du))X^{\sharp}(s,t,x)ds\\
&&-\int_{\tau}^T\nabla_y\psi(s,X(s,t,x),Y^n(s,t,x),\int_{\cO}Z^n(s,t,x,u)l(u)\nu(du))(Y^n(s,t,x))^{\sharp}ds\\
&&-\int_{\tau}^T\nabla_z\psi(s,X(s,t,x),Y^n(s,t,x),\int_{\cO}Z^n(s,t,x,u)l(u)\nu(du))\\
&&\int_{\cO}(Z^n(s,t,x,u))^{\sharp}l(u)\nu(du)ds\\
&&+\nabla\phi(X(T,t,x))(X(T,t,x))^{\sharp}.
\de

{\bf Step 4:}

By Jensen's inequality we have
\ce
&&(\int_{\cO}|Z^n(s,t,x,u)l(u)|\nu(du))^2\\
&=&(\cM(\cO))^2(\int_{\cO}|\frac{Z^n(s,t,x,u)}{l(u)}|\frac{l^2(u)\nu(du)}{\cM(\cO)})^2\\
&\leq&\cM(\cO)\int_{\cO}|Z^n(s,t,x,u)|^2\nu(du),
\de
where
\ce
\cM(\cO)=\int_{\cO}l^2(u)\nu(du).
\de
Similarly, we have
\ce
(\int_{\cO}(Z^n(s,t,x,u))^{\sharp}l(u)\nu(du))^2\leq\cM(\cO)\int_{\cO}|(Z^n(s,t,x,u))^{\sharp}|^2\nu(du).
\de
In view of Theorem \ref{BSDE existence} (2) and the assumption {\bf (BD)}, we obtain
\ce
&&\mE\widehat{\mE}[\sup_{\tau\in[0,T]}|(Y^{n+1}(\tau,t,x))^{\sharp}|^2+\int_0^T\int_{\cO}|(Z^{n+1}(\tau,t,x,u))^{\sharp}|^2\nu(du)d\tau]\\
&\leq&C\mE\widehat{\mE}[\int_0^T|\nabla_x\psi(s,X(s,t,x),Y^n(s,t,x),\int_{\cO}Z^n(s,t,x,u)l(u)\nu(du))X^{\sharp}(s,t,x)|^2ds]\\
&&+C\mE\widehat{\mE}[\int_0^T|\nabla_y\psi(s,X(s,t,x),Y^n(s,t,x),\int_{\cO}Z^n(s,t,x,u)l(u)\nu(du))(Y^n(s,t,x))^{\sharp}|^2ds]\\
&&+C\mE\widehat{\mE}[\int_0^T|\nabla_z\psi(s,X(s,t,x),Y^n(s,t,x),\int_{\cO}Z^n(s,t,x,u)l(u)\nu(du))\\
&&\int_{\cO}(Z^n(s,t,x,u))^{\sharp}l(u)\nu(du)|^2ds]\\
&&+C\mE\widehat{\mE}[|\nabla\phi(X(T,t,x))(X(T,t,x))^{\sharp}|^2]\\
&\leq&C\mE\widehat{\mE}[\int_0^T(1+\int_{\cO}|Z^n(s,t,x,u)l(u)|\nu(du))^2\\
&&(1+|X(s,t,x)|+|Y^n(s,t,x)|)^{2m}|X^{\sharp}(s,t,x)|^2ds]\\
&&+C\mE\widehat{\mE}[\int_0^T|(Y^n(s,t,x))^{\sharp}|^2ds]\\
&&+C\mE\widehat{\mE}[\int_0^T|\int_{\cO}(Z^n(s,t,x,u))^{\sharp}l(u)\nu(du)|^2ds]\\
&&+C\mE\widehat{\mE}[|(X(T,t,x))^{\sharp}|^2]\\
&\leq&C\mE\widehat{\mE}[\int_0^T(1+\int_{\cO}|Z^n(s,t,x,u)l(u)|\nu(du))^2\\
&&(1+|X(s,t,x)|+|Y^n(s,t,x)|)^{2m}|X^{\sharp}(s,t,x)|^2ds]\\
&&+C\mE\widehat{\mE}[\int_0^T|(Y^n(s,t,x))^{\sharp}|^2ds]\\
&&+C\mE\widehat{\mE}[\int_0^T\int_{\cO}|(Z^n(s,t,x,u))^{\sharp}|^2\nu(du)ds]\\
&&+C\mE\widehat{\mE}[|(X(T,t,x))^{\sharp}|^2].
\de
Now it is sufficient to deal with the term
\ce
\mE\widehat{\mE}[\int_0^T(1+\int_{\cO}|Z^n(s,t,x,u)l(u)|\nu(du))^2(1+|X(s,t,x)|+|Y^n(s,t,x)|)^{2m}|X^{\sharp}(s,t,x)|^2ds].
\de
Using H\"{o}lder's inequality, we have
\ce
&&\mE\widehat{\mE}[\int_0^T(1+\int_{\cO}|Z^n(s,t,x,u)l(u)|\nu(du))^2(1+|X(s,t,x)|+|Y^n(s,t,x)|)^{2m}|X^{\sharp}(s,t,x)|^2ds]\\
&\leq&C\mE\widehat{\mE}[\sup_{s\in[0,T]}((1+|X(s,t,x)|+|Y^n(s,t,x)|)^{2m}|X^{\sharp}(s,t,x)|^2)\\
&&(1+\int_0^T\int_{\cO}|Z^n(s,t,x,u)|^2\nu(du)ds)]\\
&\leq&C\mE\widehat{\mE}[\sup_{s\in[0,T]}((1+|X(s,t,x)|+|Y^n(s,t,x)|)^{2mp}|X^{\sharp}(s,t,x)|^{2p})]^{1/p}\\
&&\mE\widehat{\mE}[(1+\int_0^T\int_{\cO}|Z^n(s,t,x,u)|^2\nu(du)ds)^{p/(p-1)}]^{(p-1)/p}.
\de
This in turn yields the desired result by Theorem \ref{FSDE existence}.
\end{proof}

\bl\label{Picard Malliavin convergence}
Suppose $\{\widetilde{Y}(\tau,t,x),\widetilde{Z}(\tau,t,x,u);\tau\in[0,T],x\in\mR^d,u\in\cO\}$ satisfies the following backward SDEs with jumps:
\ce
&&\widetilde{Y}(\tau,t,x)+\int_{\tau}^T\int_{\cO}\widetilde{Z}(s-,t,x,u)\widetilde{N}(ds,du)\\
&=&-\int_{\tau}^T(\nabla_x\psi(s,X(s,t,x),Y(s,t,x),\int_{\cO}Z(s,t,x,u)l(u)\nu(du))X^{\sharp}(s,t,x)\\
&&+\nabla_y\psi(s,X(s,t,x),Y(s,t,x),\int_{\cO}Z(s,t,x,u)l(u)\nu(du))\widetilde{Y}(s,t,x)\\
&&+\nabla_z\psi(s,X(s,t,x),Y(s,t,x),\int_{\cO}Z(s,t,x,u)l(u)\nu(du))\int_{\cO}\widetilde{Z}(s,t,x,u)l(u)\nu(du)))ds\\
&&+\nabla\phi(X(T,t,x))X^{\sharp}(T,t,x).
\de
Then under the assumptions of Theorem \ref{BSDE existence}, we have
\ce
\mE\widehat{\mE}[\sup_{\tau\in[0,T]}|\widetilde{Y}(\tau,t,x)|^2]+\mE\widehat{\mE}[\int_0^T\int_{\cO}|\widetilde{Z}(\tau,t,x,u)|^2\nu(du)d\tau]<\infty
\de
and
\ce
&&\lim_{n\rightarrow\infty}(\mE\widehat{\mE}[\sup_{\tau\in[0,T]}|\widetilde{Y}(\tau,t,x)-(Y^{n+1}(\tau,t,x))^{\sharp}|^2]\\
&&+\mE\widehat{\mE}[\int_0^T\int_{\cO}|\widetilde{Z}(\tau,t,x,u)-(Z^{n+1}(\tau,t,x,u))^{\sharp}|^2\nu(du)d\tau])=0.
\de
\el

\begin{proof}
Using the same argument as in the proof of Lemma \ref{Picard Malliavin}, it follows that
\ce
\mE\widehat{\mE}[\sup_{\tau\in[0,T]}|\widetilde{Y}(\tau,t,x)|^2]+\mE\widehat{\mE}[\int_0^T\int_{\cO}|\widetilde{Z}(\tau,t,x,u)|^2\nu(du)d\tau]<\infty.
\de

Now we deal with the convergence of
\ce
\{(Y^{n+1}(\tau,t,x))^{\sharp},(Z^{n+1}(\tau,t,x,u))^{\sharp};\tau\in[0,T],x\in\mR^d,u\in\cO\}
\de
in $\cS^2\times\cM^2$. It is obvious that
\ce
&&\widetilde{Y}(\tau,t,x)-(Y^{n+1}(\tau,t,x))^{\sharp}+\int_{\tau}^T\int_{\cO}(\widetilde{Z}(s-,t,x,u)-(Z^{n+1}(s-,t,x,u))^{\sharp})\widetilde{N}(ds,du)\\
&=&-\int_{\tau}^T(\nabla_x\psi(s,X(s,t,x),Y(s,t,x),\int_{\cO}Z(s,t,x,u)l(u)\nu(du))X^{\sharp}(s,t,x)\\
&&+\nabla_y\psi(s,X(s,t,x),Y(s,t,x),\int_{\cO}Z(s,t,x,u)l(u)\nu(du))\widetilde{Y}(s,t,x)\\
&&+\nabla_z\psi(s,X(s,t,x),Y(s,t,x),\int_{\cO}Z(s,t,x,u)l(u)\nu(du))\int_{\cO}\widetilde{Z}(s,t,x,u)l(u)\nu(du))ds\\
&&+\int_{\tau}^T(\nabla_x\psi(s,X(s,t,x),Y^n(s,t,x),\int_{\cO}Z^n(s,t,x,u)l(u)\nu(du))X^{\sharp}(s,t,x)\\
&&+\nabla_y\psi(s,X(s,t,x),Y^n(s,t,x),\int_{\cO}Z^n(s,t,x,u)l(u)\nu(du))(Y^n(s,t,x))^{\sharp}\\
&&+\nabla_z\psi(s,X(s,t,x),Y^n(s,t,x),\int_{\cO}Z^n(s,t,x,u)l(u)\nu(du))\int_{\cO}(Z^n(s,t,x,u))^{\sharp}l(u)\nu(du))ds.
\de
Then using Theorem \ref{BSDE existence} (2), it follows that
\ce
&&\mE\widehat{\mE}[\sup_{\tau\in[0,T]}e^{\rho\tau}|\widetilde{Y}(\tau,t,x)-(Y^{n+1}(\tau,t,x))^{\sharp}|^2\\
&&+\int_0^T\int_{\cO}e^{\rho\tau}|\widetilde{Z}(\tau,t,x,u)-(Z^{n+1}(\tau,t,x,u))^{\sharp}|^2\nu(du)d\tau]\\
&\leq&\frac{C}{\rho}\mE\widehat{\mE}[\int_0^Te^{\rho s}|\nabla_x\psi(s,X(s,t,x),Y(s,t,x),\int_{\cO}Z(s,t,x,u)l(u)\nu(du))X^{\sharp}(s,t,x)\\
&&-\nabla_x\psi(s,X(s,t,x),Y^n(s,t,x),\int_{\cO}Z^n(s,t,x,u)l(u)\nu(du))X^{\sharp}(s,t,x)|^2ds]\\
&&+\frac{C}{\rho}\mE\widehat{\mE}[\int_0^Te^{\rho s}|\nabla_y\psi(s,X(s,t,x),Y(s,t,x),\int_{\cO}Z(s,t,x,u)l(u)\nu(du))\widetilde{Y}(s,t,x)\\
&&-\nabla_y\psi(s,X(s,t,x),Y^n(s,t,x),\int_{\cO}Z^n(s,t,x,u)l(u)\nu(du))(Y^n(s,t,x))^{\sharp}|^2ds]\\
&&+\frac{C}{\rho}\mE\widehat{\mE}[\int_0^Te^{\rho s}|\nabla_z\psi(s,X(s,t,x),Y(s,t,x),\int_{\cO}Z(s,t,x,u)l(u)\nu(du))\int_{\cO}\widetilde{Z}(s,t,x,u)l(u)\nu(du)\\
&&-\nabla_z\psi(s,X(s,t,x),Y^n(s,t,x),\int_{\cO}Z^n(s,t,x,u)l(u)\nu(du))\int_{\cO}(Z^n(s,t,x,u))^{\sharp}l(u)\nu(du)|^2ds]\\
&\leq&\frac{C}{\rho}\mE\widehat{\mE}[\int_0^Te^{\rho s}|\nabla_x\psi(s,X(s,t,x),Y(s,t,x),\int_{\cO}Z(s,t,x,u)l(u)\nu(du))X^{\sharp}(s,t,x)\\
&&-\nabla_x\psi(s,X(s,t,x),Y^n(s,t,x),\int_{\cO}Z^n(s,t,x,u)l(u)\nu(du))X^{\sharp}(s,t,x)|^2ds]\\
&&+\frac{C}{\rho}\mE\widehat{\mE}[\int_0^Te^{\rho s}|\nabla_y\psi(s,X(s,t,x),Y(s,t,x),\int_{\cO}Z(s,t,x,u)l(u)\nu(du))\widetilde{Y}(s,t,x)\\
&&-\nabla_y\psi(s,X(s,t,x),Y^n(s,t,x),\int_{\cO}Z^n(s,t,x,u)l(u)\nu(du))\widetilde{Y}(s,t,x)|^2ds]\\
&&+\frac{C}{\rho}\mE\widehat{\mE}[\int_0^Te^{\rho s}|\nabla_y\psi(s,X(s,t,x),Y^n(s,t,x),\int_{\cO}Z^n(s,t,x,u)l(u)\nu(du))\widetilde{Y}(s,t,x)\\
&&-\nabla_y\psi(s,X(s,t,x),Y^n(s,t,x),\int_{\cO}Z^n(s,t,x,u)l(u)\nu(du))(Y^n(s,t,x))^{\sharp}|^2ds]\\
&&+\frac{C}{\rho}\mE\widehat{\mE}[\int_0^Te^{\rho s}|\nabla_z\psi(s,X(s,t,x),Y(s,t,x),\int_{\cO}Z(s,t,x,u)l(u)\nu(du))\int_{\cO}\widetilde{Z}(s,t,x,u)l(u)\nu(du)\\
&&-\nabla_z\psi(s,X(s,t,x),Y^n(s,t,x),\int_{\cO}Z^n(s,t,x,u)l(u)\nu(du))\int_{\cO}\widetilde{Z}(s,t,x,u)l(u)\nu(du)|^2ds]\\
&&+\frac{C}{\rho}\mE\widehat{\mE}[\int_0^Te^{\rho s}|\nabla_z\psi(s,X(s,t,x),Y^n(s,t,x),\int_{\cO}Z^n(s,t,x,u)l(u)\nu(du))\int_{\cO}\widetilde{Z}(s,t,x,u)l(u)\nu(du)\\
&&-\nabla_z\psi(s,X(s,t,x),Y^n(s,t,x),\int_{\cO}Z^n(s,t,x,u)l(u)\nu(du))\int_{\cO}(Z^n(s,t,x,u))^{\sharp}l(u)\nu(du)|^2ds]\\
&\leq&\frac{C}{\rho}\mE\widehat{\mE}[\int_0^Te^{\rho s}|\nabla_x\psi(s,X(s,t,x),Y(s,t,x),\int_{\cO}Z(s,t,x,u)l(u)\nu(du))X^{\sharp}(s,t,x)\\
&&-\nabla_x\psi(s,X(s,t,x),Y^n(s,t,x),\int_{\cO}Z^n(s,t,x,u)l(u)\nu(du))X^{\sharp}(s,t,x)|^2ds]\\
&&+\frac{C}{\rho}\mE\widehat{\mE}[\int_0^Te^{\rho s}|\nabla_y\psi(s,X(s,t,x),Y(s,t,x),\int_{\cO}Z(s,t,x,u)l(u)\nu(du))\widetilde{Y}(s,t,x)\\
&&-\nabla_y\psi(s,X(s,t,x),Y^n(s,t,x),\int_{\cO}Z^n(s,t,x,u)l(u)\nu(du))\widetilde{Y}(s,t,x)|^2ds]\\
&&+\frac{C}{\rho}\mE\widehat{\mE}[\int_0^Te^{\rho s}|\widetilde{Y}(s,t,x)-(Y^n(s,t,x))^{\sharp}|^2ds]\\
&&+\frac{C}{\rho}\mE\widehat{\mE}[\int_0^Te^{\rho s}|\nabla_z\psi(s,X(s,t,x),Y(s,t,x),\int_{\cO}Z(s,t,x,u)l(u)\nu(du))\int_{\cO}\widetilde{Z}(s,t,x,u)l(u)\nu(du)\\
&&-\nabla_z\psi(s,X(s,t,x),Y^n(s,t,x),\int_{\cO}Z^n(s,t,x,u)l(u)\nu(du))\int_{\cO}\widetilde{Z}(s,t,x,u)l(u)\nu(du)|^2ds]\\
&&+\frac{C}{\rho}\mE\widehat{\mE}[\int_0^Te^{\rho s}|\int_{\cO}\widetilde{Z}(s,t,x,u)l(u)\nu(du)-\int_{\cO}(Z^n(s,t,x,u))^{\sharp}l(u)\nu(du)|^2ds]\\
&\leq&\frac{C}{\rho}\mE\widehat{\mE}[\int_0^Te^{\rho s}|\nabla_x\psi(s,X(s,t,x),Y(s,t,x),\int_{\cO}Z(s,t,x,u)l(u)\nu(du))X^{\sharp}(s,t,x)\\
&&-\nabla_x\psi(s,X(s,t,x),Y^n(s,t,x),\int_{\cO}Z^n(s,t,x,u)l(u)\nu(du))X^{\sharp}(s,t,x)|^2ds]\\
&&+\frac{C}{\rho}\mE\widehat{\mE}[\int_0^Te^{\rho s}|\nabla_y\psi(s,X(s,t,x),Y(s,t,x),\int_{\cO}Z(s,t,x,u)l(u)\nu(du))\widetilde{Y}(s,t,x)\\
&&-\nabla_y\psi(s,X(s,t,x),Y^n(s,t,x),\int_{\cO}Z^n(s,t,x,u)l(u)\nu(du))\widetilde{Y}(s,t,x)|^2ds]\\
&&+\frac{C}{\rho}\mE\widehat{\mE}[\int_0^Te^{\rho s}|\widetilde{Y}(s,t,x)-(Y^n(s,t,x))^{\sharp}|^2ds]\\
&&+\frac{C}{\rho}\mE\widehat{\mE}[\int_0^Te^{\rho s}|\nabla_z\psi(s,X(s,t,x),Y(s,t,x),\int_{\cO}Z(s,t,x,u)l(u)\nu(du))\int_{\cO}\widetilde{Z}(s,t,x,u)l(u)\nu(du)\\
&&-\nabla_z\psi(s,X(s,t,x),Y^n(s,t,x),\int_{\cO}Z^n(s,t,x,u)l(u)\nu(du))\int_{\cO}\widetilde{Z}(s,t,x,u)l(u)\nu(du)|^2ds]\\
&&+\frac{C}{\rho}\mE\widehat{\mE}[\int_0^T\int_{\cO}e^{\rho s}|\widetilde{Z}(s,t,x,u)l(u)\nu(du)-\int_{\cO}(Z^n(s,t,x,u))^{\sharp}|^2\nu(du)ds]\\
&=:&I_1^n+I_2^n+I_3^n+I_4^n+I_5^n,
\de
where we have used the fact that
\ce
&&\mE\widehat{\mE}[\int_0^Te^{\rho s}|\int_{\cO}\widetilde{Z}(s,t,x,u)l(u)\nu(du)-\int_{\cO}(Z^n(s,t,x,u))^{\sharp}l(u)\nu(du)|^2ds]\\
&\leq&C\mE\widehat{\mE}[\int_0^T\int_{\cO}e^{\rho s}|\widetilde{Z}(s,t,x,u)l(u)\nu(du)-\int_{\cO}(Z^n(s,t,x,u))^{\sharp}|^2\nu(du)ds].
\de
Using the assumption {\bf (BD)}, and the convergence of $(Y^n,Z^n)$ to $(Y,Z)$ as $n$ tends to infinity, it is easy to see that
\ce
I_1^n,I_2^n,I_4^n\rightarrow0,\quad\text{as}\quad n\rightarrow\infty.
\de
Then for any small $\epsilon>0$, we can find sufficiently large $N$ such that for all $n\geq N$,
\ce
&&\mE\widehat{\mE}[\sup_{\tau\in[0,T]}e^{\rho\tau}|\widetilde{Y}(\tau,t,x)-(Y^{n+1}(\tau,t,x))^{\sharp}|^2\\
&&+\int_0^T\int_{\cO}e^{\rho\tau}|\widetilde{Z}(\tau,t,x,u)-(Z^{n+1}(\tau,t,x,u))^{\sharp}|^2\nu(du)d\tau]\\
&&<\epsilon+\frac{C}{\rho}\mE\widehat{\mE}[\sup_{\tau\in[0,T]}e^{\rho\tau}|\widetilde{Y}(\tau,t,x)-(Y^n(\tau,t,x))^{\sharp}|^2\\
&&+\int_0^T\int_{\cO}e^{\rho\tau}|\widetilde{Z}(\tau,t,x,u)-(Z^n(\tau,t,x,u))^{\sharp}|^2\nu(du)d\tau].
\de
Now take $\rho$ sufficiently large such that
\ce
\frac{C}{\rho}<1.
\de
Then by recursion, for any $n\geq N$, we have
\ce
&&\mE\widehat{\mE}[\sup_{\tau\in[0,T]}e^{\rho\tau}|\widetilde{Y}(\tau,t,x)-(Y^{n+1}(\tau,t,x))^{\sharp}|^2\\
&&+\int_0^T\int_{\cO}e^{\rho\tau}\widetilde{Z}(\tau,t,x,u)-(Z^{n+1}(\tau,t,x,u))^{\sharp}|^2\nu(du)d\tau]\\
&<&\epsilon+\frac{C}{\rho}\mE\widehat{\mE}[\sup_{\tau\in[0,T]}e^{\rho\tau}|\widetilde{Y}(\tau,t,x)-(Y^n(\tau,t,x))^{\sharp}|^2\\
&&+\int_0^T\int_{\cO}e^{\rho\tau}|\widetilde{Z}(\tau,t,x,u)-(Z^n(\tau,t,x,u))^{\sharp}|^2\nu(du)d\tau]\\
&<&\frac{\epsilon}{1-C/\rho}+(\frac{C}{\rho})^{n+1-N}\mE\widehat{\mE}[\sup_{\tau\in[0,T]}e^{\rho\tau}|\widetilde{Y}(\tau,t,x)-(Y^N(\tau,t,x))^{\sharp}|^2\\
&&+\int_0^T\int_{\cO}e^{\rho\tau}|\widetilde{Z}(\tau,t,x,u)-(Z^N(\tau,t,x,u))^{\sharp}|^2\nu(du)d\tau].
\de
Then we complete the proof.
\end{proof}

\begin{proof}[Completion of the proof of Theorem \ref{BSDE existence}]

In view of Lemma \ref{Picard convergence}, Lemma \ref{Picard Malliavin} and Lemma \ref{Picard Malliavin convergence}, it follows that
\ce
\{Y^{n+1}(\tau,t,x),Z^{n+1}(\tau,t,x,u);\tau\in[0,T],x\in\mR^d,u\in\cO\}\in\cH_{\mD}\times\cH_{\mD,\nu},
\de
\ce
&&\lim_{n\rightarrow\infty}(\mE[\int_0^T|Y(\tau,t,x)-Y^{n+1}(\tau,t,x)|^2dt]\\
&&+\mE[\int_0^T\int_{\cO}|Z(\tau,t,x,u)-Z^{n+1}(\tau,t,x,u)|^2\nu(du)d\tau])=0
\de
and
\ce
&&\lim_{n\rightarrow\infty}(\mE\widehat{\mE}[\int_0^T|\widetilde{Y}(\tau,t,x)-(Y^{n+1}(\tau,t,x))^{\sharp}|^2dt]\\
&&+\mE\widehat{\mE}[\int_0^T\int_{\cO}|\widetilde{Z}(\tau,t,x,u)-(Z^{n+1}(\tau,t,x,u))^{\sharp}|^2\nu(du)d\tau])=0.
\de
Then using the fact that the Malliavin derivative is a closed operator, we have
\ce
\{Y(\tau,t,x),Z(\tau,t,x,u);\tau\in[0,T],x\in\mR^d,u\in\cO\}\in\cH_{\mD}\times\cH_{\mD,\nu},
\de
and the Malliavin derivative $\{Y^{\sharp}(\tau,t,x),Z^{\sharp}(\tau,t,x,u);\tau\in[0,T],x\in\mR^d,u\in\cO\}$ satisfies the following backward SDEs with jumps:
\ce
&&Y^{\sharp}(\tau,t,x)+\int_{\tau}^T\int_{\cO}Z^{\sharp}(s-,t,x,u)\widetilde{N}(ds,du)\\
&=&-\int_{\tau}^T(\nabla_x\psi(s,X(s,t,x),Y(s,t,x),\int_{\cO}Z(s,t,x,u)l(u)\nu(du))X^{\sharp}(s,t,x)\\
&&+\nabla_y\psi(s,X(s,t,x),Y(s,t,x),\int_{\cO}Z(s,t,x,u)l(u)\nu(du))Y^{\sharp}(s,t,x)\\
&&+\nabla_z\psi(s,X(s,t,x),Y(s,t,x),\int_{\cO}Z(s,t,x,u)l(u)\nu(du))\int_{\cO}Z^{\sharp}(s,t,x,u)l(u)\nu(du))ds\\
&&+\nabla\phi(X(T,t,x))X^{\sharp}(T,t,x).
\de
Therefore the proof is complete.
\end{proof}

\section{Proof of Proposition \ref{weak condition}}\label{appB}
First of all, using the condition (\ref{inverse}), there exists a constant $C_1>0$ such that
\ce
\int_{\cO}(1-e^{-\lambda|u|^3})\nu(du)\leq C_1\lambda^{\beta/3},\quad\forall\lambda>0.
\de
This implies by the inequality $xe^{-x}\leq 1-e^{-x}$ provided $x\geq0$ that
\ce
\int_{\cO}\lambda|u|^3e^{-\lambda|u|^3}\nu(du)\leq C_1\lambda^{\beta/3}.
\de
Hence for any $n\in\mN^*$, we have
\be\label{estimate 1}
\int_{\{2^{-n-1}<|u|\leq 2^{-n}\}}\lambda|u|^3e^{-\lambda|u|^3}\nu(du)\leq C_1\lambda^{\beta/3}.
\ee
On the other hand, choose $\lambda_n=2^{3n}$, then $\lambda_n|u|^3\leq1$ for any $2^{-n-1}<u\leq2^{-n}$. It follows that
\be\label{estimate 2}
\int_{\{2^{-n-1}<|u|\leq 2^{-n}\}}\lambda_n|u|^3e^{-\lambda_n|u|^3}\nu(du)\geq e^{-1}\int_{\{2^{-n-1}<|u|\leq 2^{-n}\}}2^{3n}|u|^3\nu(du).
\ee
Then using (\ref{estimate 1}) and (\ref{estimate 2}), we have
\ce
&&\int_{\{2^{-n-1}<|u|\leq 2^{-n}\}}|u|^2\nu(du)\\
&\leq&2^{n+1}\int_{\{2^{-n-1}<|u|\leq 2^{-n}\}}|u|^3\nu(du)\\
&\leq&e2^{n+1-3n}\int_{\{2^{-n-1}<|u|\leq 2^{-n}\}}\lambda_n|u|^3e^{-\lambda_n|u|^3}\nu(du)\\
&\leq&2eC_12^{n\beta-2n}=:C_22^{-n(2-\beta)}.
\de
Therefore
\ce
\int_{\{|u|\leq 2^{-n}\}}|u|^2\nu(du)=\sum_{k=n}^{\infty}\int_{\{2^{-k-1}<|u|\leq 2^{-k}\}}|u|^2\nu(du)\leq C_2\sum_{k=n}^{\infty}2^{-k(2-\beta)}\leq C_32^{-n(2-\beta)}.
\de
Now for any $\varepsilon>0$, choose $n_{\varepsilon}\in\mN^*$ such that $2^{-n_{\varepsilon}-1}<\varepsilon\leq 2^{-n_{\varepsilon}}$. This yields that
\ce
\int_{\{|u|\leq\varepsilon\}}|u|^2\nu(du)\leq\int_{|u|\leq 2^{-n_{\varepsilon}}}|u|^2\nu(du)\leq C_4\varepsilon^{2-\beta},
\de
and hence for any $p\geq2$,
\ce
\int_{\{|u|\leq\varepsilon\}}|u|^p\nu(du)\leq\varepsilon^{p-2}\int_{\{|u|\leq\varepsilon\}}|u|^2\nu(du)\leq C_4\varepsilon^{p-\beta}.
\de
Then the first conclusion of Proposition \ref{weak condition} has been proven.

Now we show that there exists a constant $C_L>0$ such that
\be\label{order condition}
\nu(|u|>\varepsilon)\leq C_L\varepsilon^{-\beta},\quad\forall \varepsilon>0.
\ee
Indeed, for any $\varepsilon>0$, choose $n_{\varepsilon}\in\mN^*$ such that $2^{-n_{\varepsilon}-1}<\varepsilon\leq2^{-n_{\varepsilon}}$. 
By the first conclusion of Proposition \ref{weak condition} we have that for any $k\in\mN^*$,
\ce
\int_{2^{-k-1}<|u|\leq 2^{-k}}|u|^2\nu(du)\leq C_O2^{-k(2-\beta)}.
\de
It follows that
\ce
\nu(|u|>\epsilon)&\leq&\sum_{k=0}^{n_{\varepsilon}}\nu(2^{-k-1}<|u|\leq 2^{-k})\\
&\leq&\sum_{k=0}^{n_{\varepsilon}}2^{2k+2}\int_{\{2^{-k-1}<|u|\leq 2^{-k}\}}|u|^2\nu(du)\\
&\leq&C_O\sum_{k=0}^{n_{\varepsilon}}2^{2k+2}2^{-k(2-\beta)}\\
&=&C_O\sum_{k=0}^{n_{\varepsilon}}2^{\beta k+2}\\
&\leq&C_52^{\beta(n_{\varepsilon}+1)}\leq C_6\varepsilon^{-\beta},
\de
and this is the desired result.

We are in a position to prove
\ce
\mE[(\int_t^{\tau}\int_{\{|u|\leq\varepsilon\}}|u|^3N(ds,du))^{-p}]\leq C_I(((\tau-t)\varepsilon^{3-\beta})^{-p}+(\tau-t)^{-\frac{3p}{\beta}}).
\de
It is easy to see that
\ce
&&\mE[e^{-\lambda\int_t^{\tau}\int_{\{|u|\leq\varepsilon\}}|u|^3N(ds,du)}]\\
&=&\exp\{-(\tau-t)\int_{\{|u|\leq\varepsilon\}}(1-e^{-\lambda|u|^3})\nu(du)\}.
\de
Then we have
\ce
&&\mE[(\int_t^{\tau}\int_{\{|u|\leq\varepsilon\}}|u|^3N(ds,du))^{-p}]\\
&=&\frac{1}{\Gamma(p)}\int_0^{\infty}\lambda^{p-1}\mE[e^{-\lambda\int_t^{\tau}\int_{\{|u|\leq\varepsilon\}}|u|^3N(ds,du)}]d\lambda\\
&=&\frac{1}{\Gamma(p)}\int_0^{\infty}\lambda^{p-1}\exp\{-(\tau-t)\int_{\{|u|\leq\varepsilon\}}(1-e^{-\lambda|u|^3})\nu(du)\}d\lambda.
\de
Applying the condition (\ref{inverse}), there exists a constant $\lambda_0>0$ such that
\ce
\int_{\cO}(1-e^{-\lambda|u|^3})\nu(du)\geq\frac{\lambda^{\beta/3}\kappa}{2},\quad\forall\lambda>\lambda_0.
\de
Now using the fact that
\ce
&&(\tau-t)\int_{\cO}(1-e^{\lambda|u|^3})\nu(du)\\
&\leq&(\tau-t)\int_{\{|u|\leq\varepsilon\}}(1-e^{-\lambda|u|^3})\nu(du)+(\tau-t)\int_{\{|u|>\varepsilon\}}(1-e^{-\lambda|u|^3})\nu(du)\\
&\leq&(\tau-t)\int_{\{|u|\leq\varepsilon\}}(1-e^{-\lambda|u|^3})\nu(du)+C_L(\tau-t)\varepsilon^{-\beta},
\de
where we have used the estimate (\ref{order condition}), for all $\lambda\geq\lambda_{\varepsilon}=(\frac{4C_L}{\kappa}\varepsilon^{-\beta})^{3/\beta}\vee\lambda_0$, we have
\ce
&&\frac{\tau-t}{\lambda^{\beta/3}}\int_{\{|u|\leq\varepsilon\}}(1-e^{-\lambda|u|^3})\nu(du)\\
&\geq&\frac{\kappa(\tau-t)}{2}-\frac{C_L(\tau-t)}{\lambda^{\beta/3}}\varepsilon^{-\beta}\geq\frac{\kappa(\tau-t)}{4}.
\de
This implies that
\ce
&&\int_{\lambda_{\varepsilon}}^{\infty}\lambda^{p-1}\exp\{-(\tau-t)\int_{\{|u|\leq\varepsilon\}}(1-e^{-\lambda|u|^3})\nu(du)\}d\lambda\\
&\leq&\int_0^{\infty}\lambda^{p-1}e^{-\frac{\kappa(\tau-t)}{4}\lambda^{\beta/3}}d\lambda\\
&\leq&C_7(\tau-t)^{-\frac{3p}{\beta}}.
\de
On the other hand, using the inequality $xe^{-x}\leq 1-e^{-x}$ provided $x\geq0$ again, and the fact that $\lambda_{\varepsilon}\varepsilon^3=\frac{4C_L}{\kappa}\vee(\lambda_0\varepsilon^3)\leq\frac{4C_L}{\kappa}\vee\lambda_0=:C_8$, we have
\ce
&&\int_0^{\lambda_{\varepsilon}}\lambda^{p-1}\exp\{-(\tau-t)\int_{\{|u|\leq\varepsilon\}}(1-e^{-\lambda|u|^3})\nu(du)\}d\lambda\\
&\leq&\int_0^{\lambda_{\varepsilon}}\lambda^{p-1}\exp\{-(\tau-t)\lambda\int_{\{|u|\leq\varepsilon\}}|u|^3 e^{-\lambda|u|^3}\nu(du)\}d\lambda\\
&\leq&\int_0^{\lambda_{\varepsilon}}\lambda^{p-1}\exp\{-(\tau-t)\lambda\int_{\{|u|\leq\varepsilon\}}|u|^3 e^{-\lambda_{\varepsilon}\varepsilon^3}\nu(du)\}d\lambda\\
&=&\int_0^{\lambda_{\varepsilon}}\lambda^{p-1}\exp\{-C_9(\tau-t)\lambda\int_{\{|u|\leq\varepsilon\}}|u|^3\nu(du)\}d\lambda\\
&\leq&\int_0^{\infty}\lambda^{p-1}\exp\{-C_9(\tau-t)\lambda\int_{\{|u|\leq\varepsilon\}}|u|^3\nu(du)\}d\lambda\\
&\leq&C_{10}((\tau-t)\int_{\{|u|\leq\varepsilon\}}|u|^3\nu(du))^{-p}\\
&\leq&C_{10}(C_O(\tau-t)\varepsilon^{3-\beta})^{-p}.
\de
Then the proof is complete.
\end{appendix}

\end{document}